\documentclass{amsart}
\usepackage{amsmath,amssymb,bm}
\usepackage{amsmath}
\usepackage{amssymb}
\usepackage{amscd}
\usepackage{amsthm}
\usepackage[dvips]{graphicx}

\theoremstyle{definition}
\newtheorem{defn}{Definition}[section]
\newtheorem{example}[defn]{Example}
\newtheorem{rem}[defn]{Remark}
\theoremstyle{plain}
\newtheorem{thm}[defn]{Theorem}
\newtheorem{prop}[defn]{Proposition}
\newtheorem{lem}[defn]{Lemma}

\newcommand{\Aut}{\operatorname{Aut}}
\newcommand{\Hom}{\operatorname{Hom}}
\newcommand{\id}{\operatorname{id}}

\numberwithin{equation}{section}
\title{Unoriented HQFT and its underlying algebra}
\author{Keiji Tagami}
\date{\today}
\address{
Department of Mathematics,
Tokyo Institute of Technology,
Oh-okayama, Meguro, Tokyo 152-8551, Japan
}
\email{tagami.k.aa@m.titech.ac.jp}
\begin{document}
\maketitle
\begin{abstract}
Turaev and Turner introduced a bijection between unoriented topological quantum field theories and extended Frobenius algebras. 
In this paper, we will show that there exists a bijective correspondence between unoriented ($1+1$)-dimensional homotopy quantum field theories and extended crossed group algebras. 
\end{abstract}
\section{Introduction}
\par
In \cite{atiyah1}, Atiyah introduced a mathematical definition of topological quantum field theories (TQFTs). 
A ($d+1)$-TQFT assigns a module to each $d$-dimensional manifold and assigns a homomorphism of modules to each ($d+1$)-dimensional cobordism. 
Abrams \cite{abrams:1996} showed that there is a bijective correspondence between oriented (1+1)-TQFTs and Frobenius algebras. 
Turaev \cite{turaev:1999} defined a concept of homotopy quantum field theories (HQFTs) with target $X$, where $X$ is a connected topological space. 
An HQFT assigns a module and a homomorphism of modules to each ``$X$-manifold'' and ``$X$-cobordism" respectively. 
For any group $\pi$, he constructed a bijective correspondence between oriented ($1+1$)-dimensional HQFTs with target $X$ for $X=K(\pi, 1)$ and crossed $\pi$ algebras in \cite{turaev:1999}, where a crossed $\pi$ algebra $V$ is a Frobenius $\pi$-algebra endowed with a group homomorphism  $\varphi\colon\pi\rightarrow  \Aut (V)$. 
In \cite{staic-turaev:2009} Staic and Turaev discussed ($1+1$)-dimensional HQFTs more generally. 
Turaev and Turner \cite{turner-turaev:2006} showed that there exists a bijective correspondence between unoriented (1+1)-TQFTs and extended Frobenius algebras. 
An extended Frobenius algebra $K$ is a Frobenius algebra endowed with an element $\theta \in K$ and a homomorphism $\Phi\colon K \rightarrow K$. 
\par
In this paper, we consider a group $\pi$ such that $\alpha^{2}=1$ for any $\alpha\in\pi$, $X=K(\pi, 1)$ and unoriented (1+1)-dimensional HQFTs with target $X$. 
Note that such a group $\pi$ is a $\mathbf{Z}/2\mathbf{Z}$ vector space. 
Moreover we introduce ``extended crossed $\pi$-algebra'' $L$ which consists of a Frobenius $\pi$-algebra, a group homomorphism  $\varphi\colon\pi\rightarrow \Aut(L)$, elements $\{\theta_{\alpha} \in L\}_{\alpha\in\pi}$ and a homomorphism $\Phi \colon L\rightarrow L$ (Definition~$\ref{extcross}$). 
We will show that there is a bijective correspondence between unoriented (1+1)-dimensional HQFTs with target $X$ and extended crossed $\pi$-algebras (Theorem~$\ref{mainthm}$). 
\par
In Section~$\ref{preliminaries}$, we recall definitions of HQFTs and some algebras introduced in \cite{turaev:1999} and will define unoriented HQFTs and extended crossed group algebras. 
In Section~$\ref{underlying}$, we construct an extended crossed group algebra from an HQFT ($A, \tau$). We call it underlying extended crossed group algebra of ($A, \tau$). At the end of this section, we introduce our main theorem  (Theorem~$\ref{mainthm}$). 
In Sections~$\ref{main}$ and $\ref{main2}$, we prove the main theorem. 
In Section~$\ref{ex}$, we give some examples. 
\par
Throughout this paper, the symbol $R$ denotes a commutative ring with unit and the symbol $\pi$ denotes a group.  

\section{Unoriented HQFTs and extended crossed group algebras}\label{preliminaries}
Here we will explain terminology used in this paper. 
\subsection{Unoriented HQFTs}
In this subsection, we recall the definition of unoriented homotopy quantum field theories. An oriented homotopy quantum field theory is introduced by Turaev \cite{turaev:1999}.  
\begin{defn}[\cite{turaev:1999}]
Let $X$ be a $K(\pi$, $1)$ space with a base point $x_{0}\in X$. 
A pair {\rm(}$M$, $g_{M}${\rm)} is called an {\it unoriented $X$-manifold} if $M$ is a pointed closed unoriented manifold and $g_{M}$ is a map from $M$ to $X$. 
We call the map $g_{M}$ the {\it characteristic map}. 
Since the spaces M and X are pointed, the map $g_{M}$ sends the base points of all components of $M$ to $x_{0}$. 
A disjoint union of unoriented $X$-manifolds and the empty set are also unoriented $X$-manifolds. 
An {\it unoriented $X$-homeomorphism} of unoriented $X$-manifolds $f\colon(M, g_{M})\rightarrow (M', g_{M'})$ is a homeomorphism from $M$ to $M'$ sending the base points of $M$ to those of $M'$ such that $g_{M}=g_{M'}\circ f$.
\end{defn}
\par
\begin{defn}[\cite{turaev:1999}]
Let $X$ be a $K(\pi$, $1)$ space with a base point $x_{0}\in X$. 
An {\it unoriented $X$-cobordism} is a tuple {\rm(}$W, M_{0}, M_{1}, g${\rm)} such that the triple {\rm(}$W, M_{0}, M_{1}${\rm)} is an unoriented cobordism, that $M_{0}$ and $M_{1}$ are unoriented $X$-manifolds and that $g\colon W\rightarrow X$ is a map which sends the base points of $M_{0}$ and $M_{1}$ to $x_{0}\in X$. 
We call the boundary $M_{0}$ the bottom base, $M_{1}$ the top base and the map $g$ the characteristic map. 
An {\it unoriented $X$-homeomorphism} of $X$-cobordisms $f\colon (W, M_{0}, M_{1}, g)\rightarrow (W', M'_{0}, M'_{1}, g')$ is a homeomorphism from $W$ to $W'$ inducing unoriented $X$-homeomorphisms $M_{0}\rightarrow M'_{0}$ and $M_{1}\rightarrow M'_{1}$ such that $g=g'\circ f$. 
\end{defn}
\begin{defn}[\cite{turaev:1999}]\label{HQFT}
Fix an integer $d \geq 0$ and a path connected topological space $X$ with base point $x \in X$. An {\it unoriented {\rm(}$d+1${\rm)}-dimensional homotopy quantum field theory} {\rm(}HQFT for short {\rm)} {\rm(}$A, \tau${\rm)} over $R$ with target $X$ assigns 
\begin{itemize}
\item a finitely generated projective $R$-module $A(M, g)$ {\rm(}$A(M)$ for short {\rm)} to any unoriented $d$-dimensional $X$-manifold $(M, g)$, 
\item an $R$-isomorphism $f_{\sharp}\colon A(M,g)\rightarrow A(M', g')$ to any unoriented $X$-homeomorphism of $d$-dimensional $X$-manifolds $f\colon (M, g)\rightarrow (M', g')$, 
\item an $R$-homomorphism $\tau(W, g)\colon A(M_{0}, g|_{M_{0}})\rightarrow A(M_{1}, g|_{M_{1}})$ to any {\rm(}$d+1${\rm)}-dimensional $X$-cobordism {\rm(}$W, M_{0}, M_{1}, g${\rm)}. 
\end{itemize}
Moreover these modules and homomorphisms should satisfy the following axioms: 
\par
$(1)$ for unoriented $X$-homeomorphisms of unoriented $X$-manifolds $f\colon M\rightarrow M'$ and $f'\colon M'\rightarrow M''$, we have $(f'\circ f)_{\sharp}=f'_{\sharp}\circ f_{\sharp}$, 
\par
$(2)$ for unoriented $d$-dimensional $X$-manifolds $M$ and $N$, there is a natural isomorphism $A(M\sqcup N)=A(M)\otimes A(N)$, where $M\sqcup N$ is the disjoint union of $M$ and $N$, 
\par
$(3)$ $A(\emptyset)=R$, 
\par
$(4)$ for any unoriented $X$-cobordism W, the homomorphism $\tau(W)$ is natural with respect to unoriented $X$-homeomorphisms, 
\par
$(5)$ if an unoriented $(d+1)$-dimensional $X$-cobordism $(W, M_{0}, M_{1}, g)$ is the disjoint union of two unoriented $(d+1)$-dimensional $X$-cobordisms $W_{0}$ and $W_{1}$, then $\tau(W)=\tau(W_{1})\otimes \tau(W_{0})$, 
\par
$(6)$ if an oriented $(d+1)$-dimensional $X$-cobordism $(W, M_{0}, M_{1}, g)$ is obtained from two $(d+1)$-dimensional $X$-cobordism $(W_{0}, M_{0}, N)$ and $(W_{1}, N', M_{1})$ by gluing along $f\colon N\rightarrow N'$, then $\tau(W)=\tau(W_{1})\circ f_{\sharp} \circ \tau(W_{0})$, 
\par
$(7)$ for any unoriented $d$-dimensional $X$-manifold {\rm(}$M, g${\rm)} and any continuous map $F\colon M\times [0,1]\rightarrow X$ such that $F|_{M\times 0}=F|_{M\times 1}=g$ and that $F(m\times [0,1])=\{x\}$ for any base point $m$ of $M$, we have $\tau(M\times [0,1], M\times 0, M\times 1, F)=\id_{A(M)}\colon A(M)\rightarrow A(M)$, 
\par
$(8)$ for any unoriented $(d+1)$-dimensional X-cobordism {\rm(}$W, g${\rm)}, $\tau(W)$ is preserved under any homotopy of $g$ relative to $\partial{W}$. 
\end{defn}
If two maps $f$ and $f'\colon M\rightarrow X$ are homotopic, there is a natural isomorphism $A(M, f)\cong A(M, f')$. Hence we can suppose that $A(M, f)$ is preserved under any homotopy of $f$. Similary $\tau(W, g)$ is preserved under any homotopy of $g$ (maybe not relative to $\partial{W}$). 
\subsection{Extended crossed group algebras}
In this subsection, we recall some algebras which are introduced in \cite{turaev:1999} and define extended crossed group algebras. 
\begin{defn}
An {\it $R$-algebra} $L$ is a $\pi$-algebra over the ring $R$ if $L$ is an associative algebra over $R$ endowed with a splitting $L=\bigoplus_{\alpha  \in \pi}L_{\alpha}$ such that each $L_{\alpha}$ is a finitely generated projective $R$-module, that $L_{\alpha}L_{\beta}\subset L_{\alpha \beta }$ for any $\alpha, \beta \in\pi$, and that $L$ has the unit element $1_{L}\in L_{1}$. 
\end{defn}
Let $V$ and $W$ be $R$-modules and $\eta\colon V\otimes W\rightarrow R$ be a bilinear form. The map $\eta$ is non-degenerate if the two maps $d\colon V\rightarrow \Hom_{R}(W, R)$ defined by $d(v)(w):=\eta(v, w)$ and $s\colon W\rightarrow \Hom_{R}(V, R)$ defined by $s(w)(v):=\eta(v, w)$ are isomorphisms, where $v\in V$ and $w\in W$.  
\begin{defn}[\cite{turaev:1999}]
A pair $(L, \eta)$ is a {\it Frobenius $\pi$-algebra} over $R$ if $L$ is a $\pi$-algebra over $R$ and $\eta \colon L_{\alpha }\otimes L_{\beta}\rightarrow R $ is an $R$-bilinear form such that 
\par
$(1)$
$\eta(L_{\alpha}\otimes L_{\beta})=0$ if $\alpha \beta \neq 1$ and the restriction of $\eta$ to $L_{\alpha }\otimes L_{\alpha^{-1} }$ is non-degenerate for any $\alpha \in\pi$, 
\par
$(2)$
$\eta(ab, c)=\eta(a, bc)$ for any $a, b, c\in L$. 
\end{defn}
A Frobenius $\pi$-algebra with $\pi$ a trivial group is called a {\it Frobenius algebra} (\cite{abrams:1996}). 
\par
For any Frobenius $\pi$-algebra $(L, \eta)$, $\Aut(L)$ is a group which consists of algebra automorphisms preserving $\eta$. 
\begin{defn}[\cite{turaev:1999}]
A triple $(L, \eta, \varphi)$ is a {\it crossed $\pi$-algebra} over $R$ if the pair $(L, \eta)$ is a Frobenius $\pi$-algebra over $R$ and $\varphi\colon \pi\rightarrow\Aut(L)$ is a group homomorphism satisfying the following axioms:
\par
$(1)$
for any $\beta\in\pi$, $\varphi_{\beta}:=\varphi(\beta)$ satisfies $\varphi_{\beta}(L_{\alpha})\subset L_{\beta \alpha \beta ^{-1}}$ for any $\alpha \in\pi$, 
\par
$(2)$
$\varphi_{\alpha}|_{L_{\alpha}}=\id_{L_{\alpha}}$ for any $\alpha \in\pi$, 
\par
$(3)$
for any $a\in L_{\alpha}$ and $b\in L_{\beta}$, we have $\varphi_{\beta}(a)b=ba$, 
\par
$(4)$
for any $\alpha, \beta\in\pi$ and any $c\in L_{\alpha \beta \alpha^{-1} \beta^{-1} }$, we have $\operatorname {Tr}(c\varphi_{\beta}\colon L_{\alpha }\rightarrow L_{\alpha })=\operatorname {Tr}(\varphi_{\alpha^{-1}}c\colon L_{\beta}\rightarrow L_{\beta})$, where {\rm Tr} is the $R$-valued trace of endmorphisms of finitely generated projective $R$-modules (see for instance \cite{quantum-turaev}). 
\end{defn}
In \cite{turaev:1999}, Turaev showed that there exists a relation between oriented HQFTs with target $K(\pi, 1)$ space and crossed $\pi$-algebras. 
\begin{thm}[Theorem~$4.1$ in \cite{turaev:1999}]\label{turaev}
Let $\pi$ be a group and $X$ be a $K(\pi, 1)$ space. Then every oriented $(1+1)$-dimensional HQFT with target $X$ over the ring $R$ determines an underlying crossed $\pi$-algebra over $R$. This induces a bijection between the set of isomorphism classes of oriented $(1+1)$-dimensional HQFTs and the set of isomorphism classes of crossed $\pi$-algebras. 
\end{thm}
For any crossed $\pi$-algebra $(L, \eta, \varphi)$, we denote the HQFT corresponding to the crossed $\pi$-algebra by $(A^{L}, \tau^{L})$ . 
Now we define extended crossed group-algebras. 

\begin{figure}[!h]
\begin{center}
\includegraphics[scale=0.18]{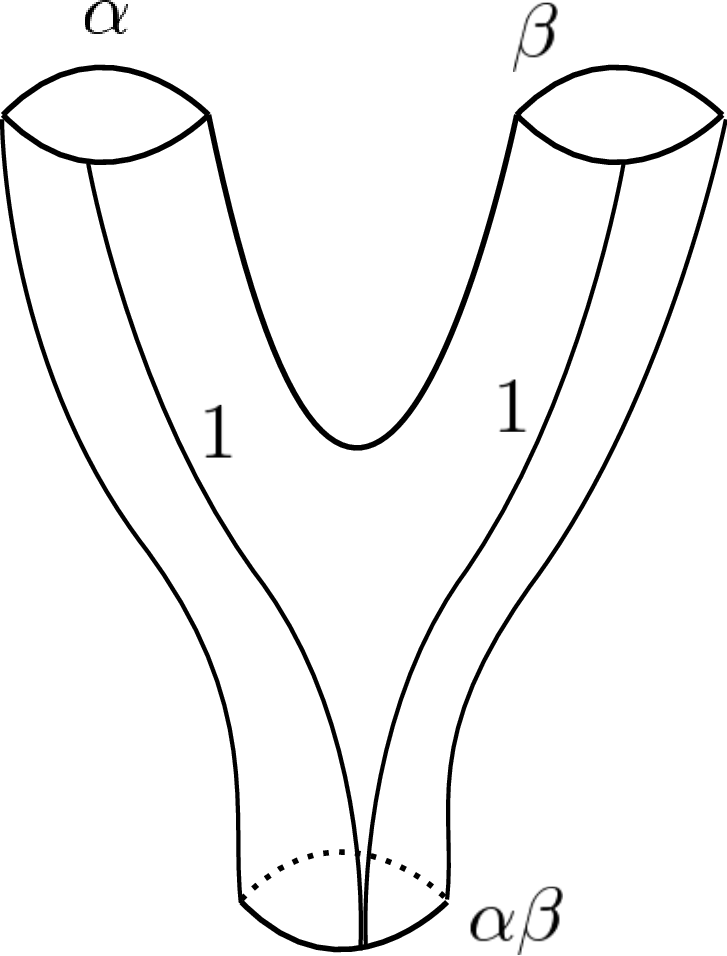}
\end{center}
\caption{The cobordism $D_{+,+,-}$($\alpha, \beta; 1, 1$). }
\label{yoseki}
\end{figure}

\begin{figure}[!h]
\begin{center}
\includegraphics[scale=0.13]{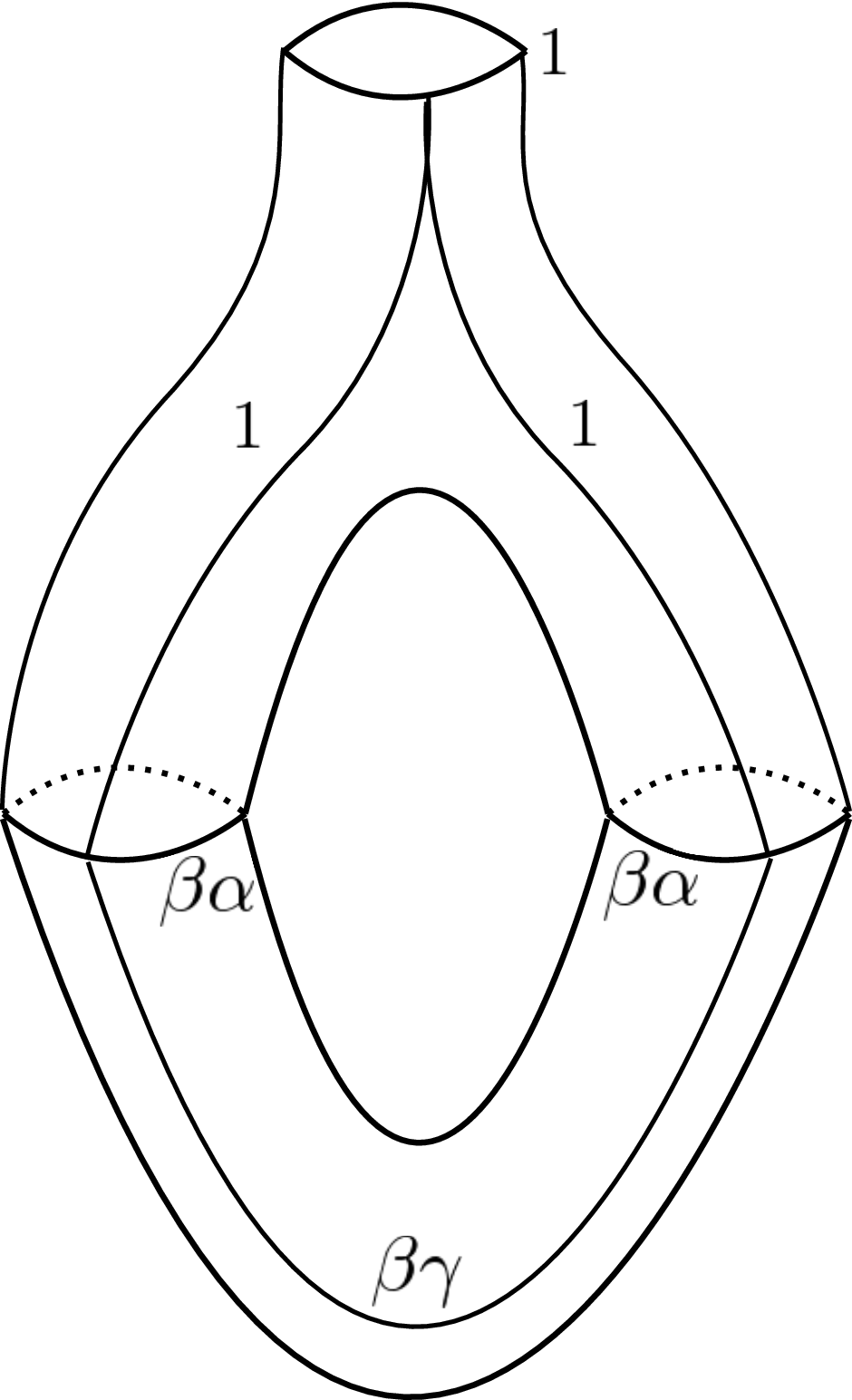}
\end{center}
\caption{Definition of the cobordism $(Q, \emptyset, (\mathbf{S}^{1}, 1))$. }
\label{tennkaizu1}
\end{figure}

\begin{figure}[!h]
\begin{center}
\includegraphics[scale=0.22]{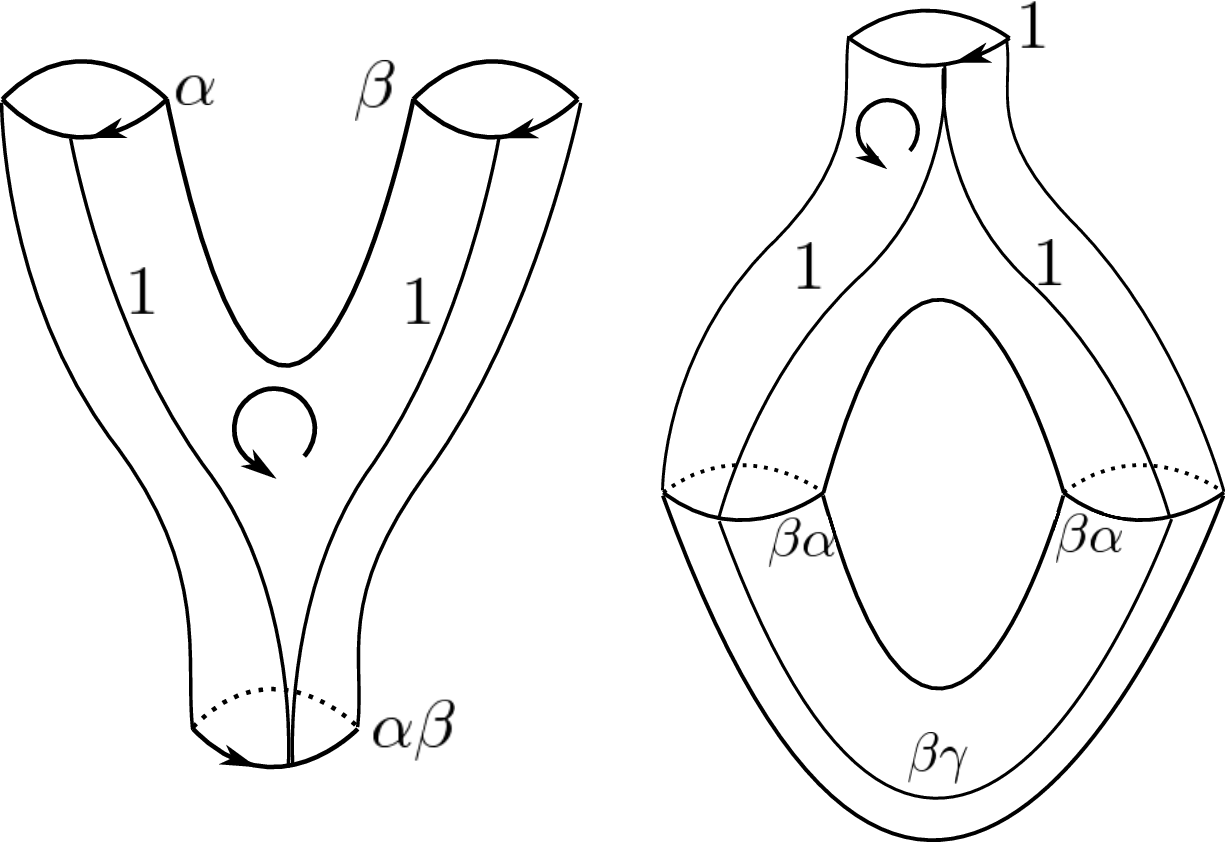}
\end{center}
\caption{Orientations.}
\label{muki}
\end{figure}

\begin{figure}[!h]
\begin{center}
\includegraphics[scale=0.3]{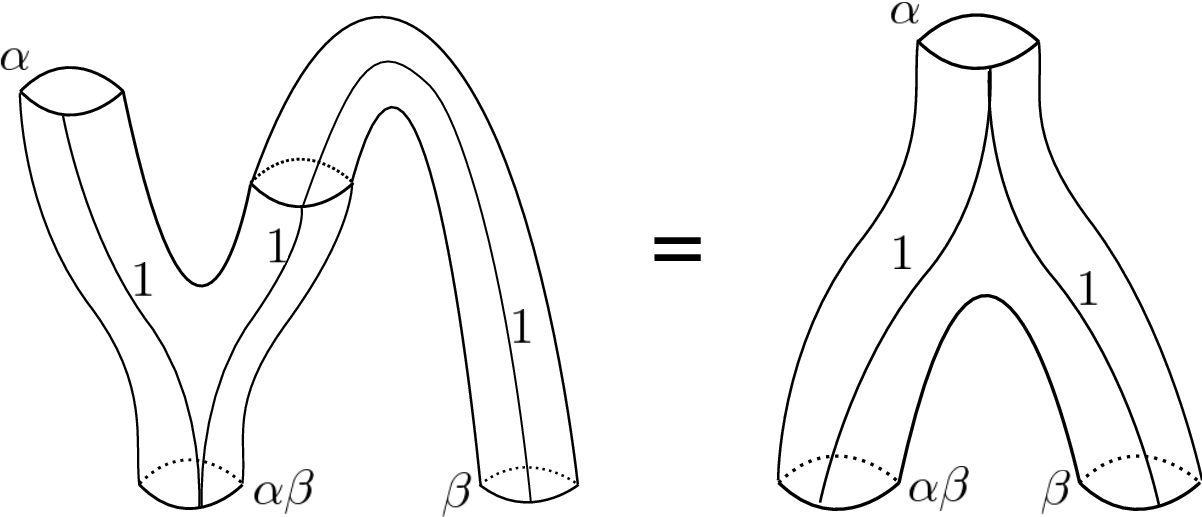}
\end{center}
\caption{Relation of cobordisms $(id\otimes\eta)\circ(\Delta _{\alpha , \beta}\otimes\id)=m$.}
\label{computeyoseki}
\end{figure}

\begin{figure}[!h]
\begin{center}
\includegraphics[scale=0.3]{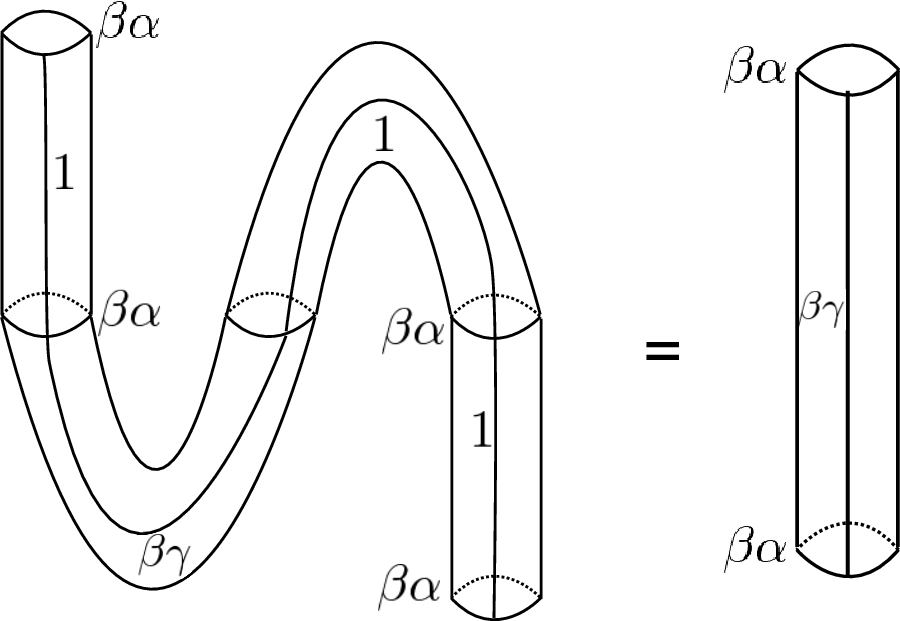}
\end{center}
\caption{Relation of cobordisms $(\id\otimes\eta)\circ\tau(Q')\otimes \id=\varphi _{\beta \gamma }$.}
\label{comp2}
\end{figure}

\begin{defn}\label{extcross}
Let $\pi$ be a group such that $\alpha ^{2}=1$ for any $\alpha \in\pi$. A tuple $(L, \eta, \varphi, \{\theta_{\alpha}\}_{\alpha\in\pi}, \Phi)$ is an {\it extended crossed $\pi$-algebra} over $R$ if the triple $(L, \eta, \varphi)$ is a crossed $\pi$-algebra, and the family of elements $\{\theta_{\alpha}\in L_{1}\}_{\alpha\in\pi}$ and the homomorphism of $R$-modules $\Phi\colon L\rightarrow L$ satisfy the following axioms:
\par
$(1)$
$\Phi^{2}=\id$, 
\par
$(2)$
$\Phi(L_{\alpha})\subset L_{\alpha}$ for any $\alpha\in\pi$, 
\par
$(3)$
for any $v, w\in L$, $\Phi(vw)=\Phi(w)\Phi(v)$, 
\par
$(4)$
$\Phi(1_{L})=1_{L}$, 
\par
$(5)$
$\eta\circ(\Phi\otimes \Phi)=\eta$, 
\par
$(6)$
for any $\alpha \in\pi$, $\Phi\circ\varphi_{\alpha}=\varphi_{\alpha}\circ\Phi$, 
\par
$(7)$
for any $\alpha, \beta, \gamma\in\pi$ and $v\in L_{\alpha \beta }$, we have 
\begin{center}
$m\circ(\Phi\otimes \varphi_{\gamma })\circ \Delta _{\alpha, \beta }(v)=\varphi_{\gamma }(\theta_{\alpha \gamma }\theta_{\gamma }v)$, 
\par
$m\circ(\varphi_{\gamma }\otimes \Phi)\circ \Delta _{\alpha, \beta }(v)=\varphi_{\gamma }(\theta_{\beta \gamma }\theta_{\gamma }v)$, 
\end{center}
where $\Delta _{\alpha, \beta }\colon L_{\alpha \beta }\rightarrow L_{\alpha}\otimes L_{\beta}$ is defined by the following relation: 
\begin{align}
(\id\otimes\eta)\circ(\Delta_{\alpha , \beta }\otimes\id)=m. \label{eq2-1}
\end{align}
Since $\eta$ is non-degenerate and each $L_{\alpha }$ is finitely generated, such a map $\Delta_{\alpha ,\beta } $ is uniquely determined. 
\par
$(8)$
for any $\alpha, \beta \in\pi$ and $v\in L_{\alpha }$, we have $\Phi(\theta_{\beta}v_{\alpha})=\varphi_{\beta\alpha}(\theta_{\beta\alpha}v_{\alpha})$, 
\par
$(9)$
for any $\alpha\in\pi$, we have $\Phi(\theta_{\alpha})=\theta_{\alpha}$, 
\par
$(10)$
for any $\alpha, \beta\in\pi$, we have $\varphi_{\beta}(\theta_{\alpha})=\theta_{\alpha}$, 
\par
$(11)$ for any $\alpha, \beta, \gamma\in\pi$, we have $\theta_{\alpha }\theta_{\beta }\theta_{\gamma }=q(1)\theta_{\alpha \beta \gamma }$, where $q\colon R\rightarrow L_{1}$ is defined as follows. 
Let $\{a_{i}\in L_{\alpha \beta }\}_{i=1}^{n}$ and $\{b_{i}\in L_{\alpha \beta }\}_{i=1}^{n}$ be families of elements of $L_{\alpha \beta }$ satisfying the 
following condition: for any $v\in L_{\alpha \beta }$ 
\begin{align}
\sum_{i}\eta(b_{i}\otimes v)a_{i}=\varphi_{\beta \gamma }(v). \label{eq2-2}
\end{align}
From the same reason as ($7$), such $a_{i}$ and $b_{i}$ are uniquely determined. Then we put $q(1):=\sum_{i}a_{i}b_{i}$.  
\end{defn}
\begin{rem}
(1)
Let $D_{+,+,-}(\alpha, \beta;1, 1)$ be the oriented $X$-cobordism given by Figure~$\ref{yoseki}$. Its bottom base is a $X$-manifold $(\mathbf{S}^{1}, \alpha \beta )$ and its top base is the disjoint union of two $X$-manifolds $(\mathbf{S}^{1}, \alpha)$ and $(\mathbf{S}^{1}, \beta )$. Its characteristic map sends each labeled arc to the loop corresponding to the label. 
Such a map is uniquely determined up to homotopy since $X$ is $K(\pi, 1)$ space. 
The orientation of $D_{+,+,-}(\alpha, \beta;1, 1)$ is given by Figure~$\ref{muki}$. Then we have $\tau^{L}(D_{+,+,-}(\alpha, \beta;1, 1))=\Delta _{\alpha, \beta}$. The relation ($\ref{eq2-1}$) corresponds to Figure~$\ref{computeyoseki}$. 
\par
(2)
Let $Q$ be the $X$-cobordism depicted in Figure $\ref{tennkaizu1}$. 
It is a once-punctured torus whose bottom base is empty and whose top base is a $X$-manifold $(\mathbf{S}^{1}, 1)$. 
Its characteristic map sends each labeled arc to the loops corresponding to the label. 
Its orientation is given by Figure~$\ref{muki}$. Then we have $q=\tau^{L}(Q)$. The relation ($\ref{eq2-2}$) corresponds to Figure~$\ref{comp2}$. 
\end{rem}

\section{Underlying algebraic structures of HQFTs}\label{underlying}
\begin{figure}[!h]
\begin{center}
\includegraphics[scale=0.2]{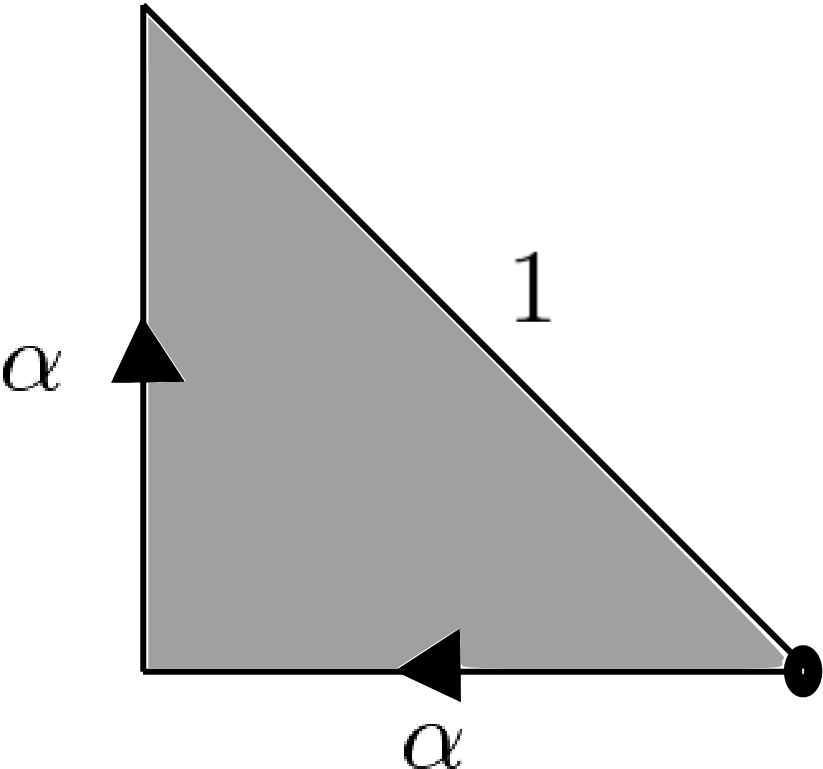}
\end{center}
\caption{The cobordism $(Mb, \emptyset, \partial (Mb), g_{\alpha})$. Three edges are labeled by $\alpha, \alpha, 1\in\pi$. 
Two edges of the triangle labeled by $\alpha\in \pi$ are identified along the arrows depicted in Figure~$\ref{Mobius band}$. 
The map $g_{\alpha}$ is a continuous map from $Mb$ to $X$ sending each edge to the loop corresponding to the label. 
}
\label{Mobius band}
\end{figure}

In this section, we construct an extended crossed group algebra from an HQFT. 
Assume that $\pi$ is a group such that any element $\alpha\in \pi$ satisfies $\alpha^{2}=1 (=1_{\pi})$, where $1_{\pi}$ is the unit of $\pi$ (in particular, $\pi$ is an abelian group). 
Moreover let $X$ be a $K(\pi, 1)$ space with a base point $x_{0}\in X$. 
Throughout this section, let $(A, \tau)$ be an unoriented ($1+1$)-dimensional HQFT with target $X$. 
Let $\mathbf{S}^{1}$ be an (unoriented) circle. 
For any unoriented $1$-dimensional $X$-manifold $(\mathbf{S}^{1}, g)$, if we give $\mathbf{S}^{1}$ an orientation, we can regard the homotopy class of $g$ as an element $\alpha\in \pi=\pi_{1}(X)$. 
The element $\alpha$ does not depend on the choice of the orientation of $\mathbf{S}^{1}$ since $\alpha=\alpha^{-1}$. 
Since we consider the module $A(\mathbf{S}^{1}, g)$, we can denote the unoriented $1$-dimensional $X$-manifold $(\mathbf{S}^{1}, g)$ by $(\mathbf{S}^{1}, \alpha)$. 
\begin{defn}
Let $Mb$ be a M\"{o}bius band. 
For any $\alpha\in \pi$, we define an unoriented $(1+1)$-dimensional $X$-cobordism $(Mb, \emptyset, \partial (Mb), g_{\alpha})$ as the unoriented $(1+1)$-dimensional $X$-cobordism in Figure~$\ref{Mobius band}$.
Choose an unoriented X-homeomorphism $f\colon (\partial(Mb), g|_{\partial(Mb)})\rightarrow (\mathbf{S}^{1}, 1)$, and define an element $\theta_{\alpha}$ by
\par
\begin{center}
$\theta_{\alpha}:=f_{\sharp}(\tau((Mb, \emptyset, \partial (Mb), g_{\alpha}))(1))\in A(\mathbf{S}^{1}, 1) $. 
\end{center}
\end{defn}
\begin{lem}\label{lem1}
The element $\theta_{\alpha}$ does not depend on the choice of an unoriented X-homeomorphism $f$. 
\end{lem}
\begin{proof}
Let $f$ and $f'\colon (\partial(Mb), g|_{\partial(Mb)})\rightarrow (\mathbf{S}^{1}, 1)$ be unoriented $X$-homeomorphisms such that $f$ is not isotopic to $f'$. 
Let $T\colon Mb\rightarrow Mb$ be a homeomorphism reversing the orientation of the boundary. 
The map $T$ induces an unoriented $X$-homeomorphism $T\colon (Mb, g_{\alpha})\rightarrow (Mb, g_{\alpha})$. 
Then $f\circ T|_{\partial(Mb)}$ is isotopic to $f'$. 
By Definition~$\ref{HQFT}$, we have 
\begin{align*}
f_{\sharp}(\tau(Mb, g_{\alpha})(1))
&=f_{\sharp}(\tau(T(Mb), g_{\alpha})(1))\\
&=f_{\sharp}\circ(T|_{\partial(Mb)})_{\sharp}(\tau(Mb, g_{\alpha})(1))\\
&=f'_{\sharp}(\tau(Mb, g_{\alpha})(1)). 
\end{align*}
Therefore the element $\theta_{\alpha}$ does not depend on the choice of $f$. 
\end{proof}
%
%
%
%
%
%
\par
Let $\chi \colon \mathbf{S}^{1}\rightarrow \mathbf{S}^{1}$ be a homeomorphism reversing the orientation. 
For any $\alpha\in\pi$, we define an isomorphism of $R$-modules $\Phi_{\alpha}=\Phi\colon A(\mathbf{S}^{1}, \alpha)\rightarrow A(\mathbf{S}^{1}, \alpha)$ by $\chi_{\sharp}:A(\mathbf{S}^{1}, \alpha)\rightarrow A(\mathbf{S}^{1}, \alpha)$, where $\chi_{\sharp}$ is the $R$-homomorphism induced by $\chi$. 
Clearly we have $\Phi(\theta_{\alpha})=\theta_{\alpha}$ for any $\alpha\in \pi$. 
%
%
%
\par
\begin{figure}[!h]
\begin{center}
\includegraphics[scale=0.15]{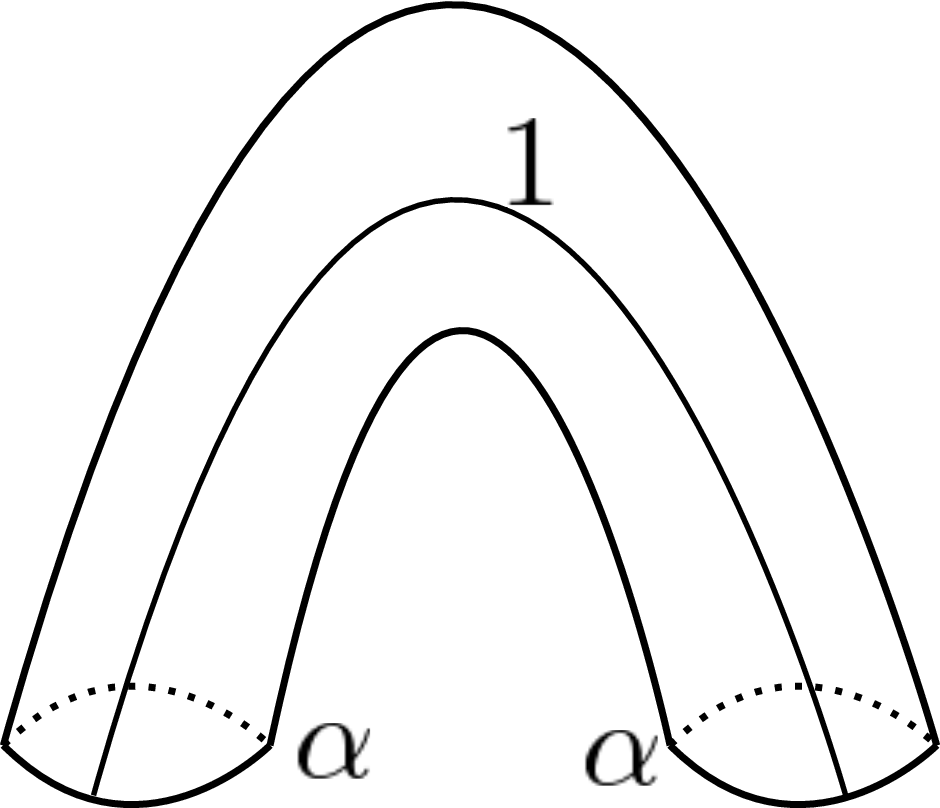}
\end{center}
\caption{The cobordism $C_{-,-}$($\alpha; 1$). }
\label{e-ta}
\end{figure}

For any $\alpha\in\pi$, let $C_{-,-}(\alpha; 1)$ be an unoriented $X$-cobordism depicted in Figure~$\ref{e-ta}$. 
The unoriented $X$-cobordism $C_{-,-}(\alpha; 1)$ is an annulus whose bottom base is the disjoint union of two copies of $X$-manifolds $(\mathbf{S}^{1}, \alpha)$, whose top base is empty and whose characteristic map sends the arc labeled by $1\in\pi$ onto $x_{0}\in X$. 
For any element $\alpha\in \pi$, we define a homomorphism of $R$-modules $\eta _{\alpha}=\eta \colon A(\mathbf{S}^{1}, \alpha)\otimes A(\mathbf{S}^{1}, \alpha)\rightarrow R$ by $\tau(C_{-,-}(\alpha, 1))$. 
\begin{lem}\label{lem3}
We have $\eta \circ(\Phi \otimes\Phi )=\eta$. 
\end{lem}
\begin{proof}
Let a map $\mu\colon C_{-,-}(\alpha; 1)\rightarrow C_{-,-}(\alpha; 1)$ be an orientation reversing homeomorphism. From Definition~$\ref{HQFT}$, we have 
\begin{center}
$\tau(C_{-,-}(\alpha, 1))=\tau(\mu (C_{-,-}(\alpha, 1)))\circ (\mu |_{(\mathbf{S}^{1}, \alpha)\sqcup(\mathbf{S}^{1}, \alpha)})_{\sharp}=\tau(C_{-,-}(\alpha, 1))\circ (\Phi \otimes\Phi )$. 
\end{center}
\end{proof}
\par
\begin{figure}[!h]
\begin{center}
\includegraphics[scale=0.15]{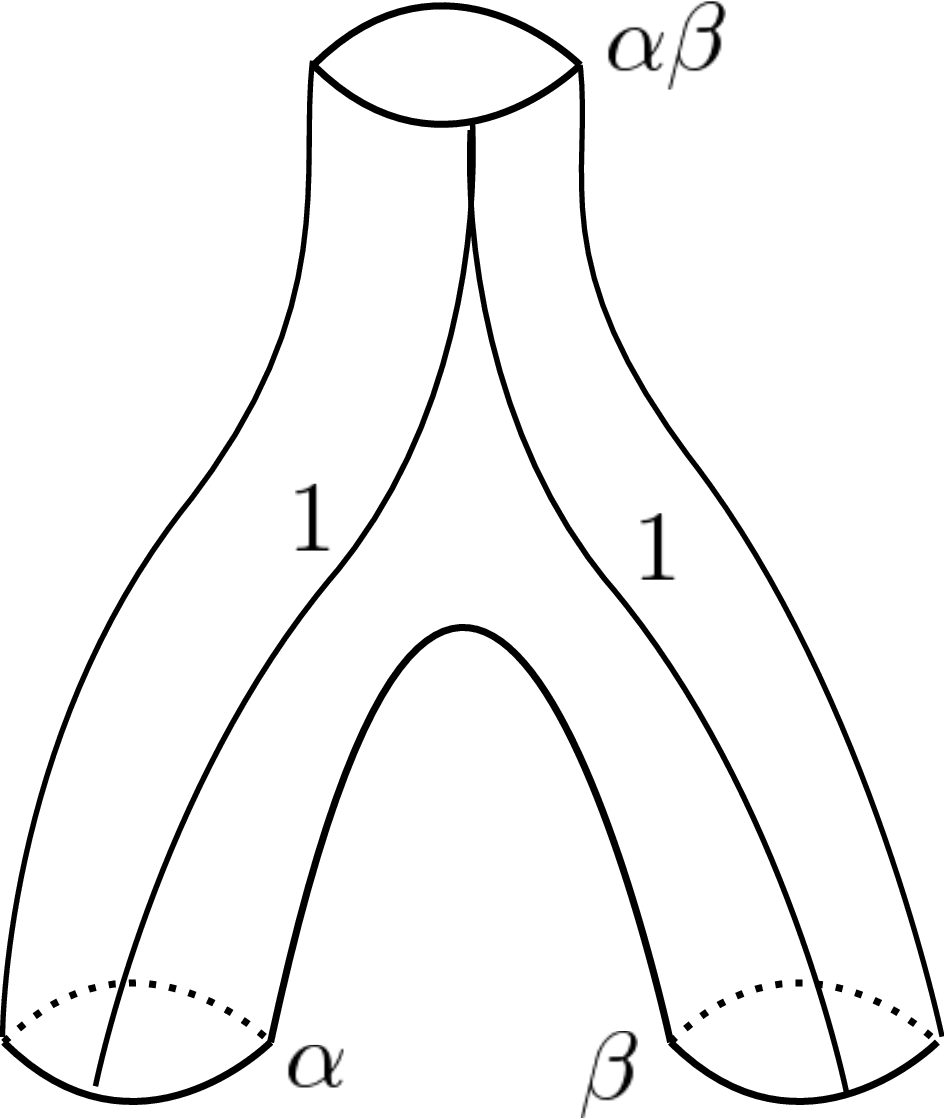}
\end{center}
\caption{The cobordism $D_{-,-,+}$($\alpha, \beta; 1, 1$).}
\label{seki}
\end{figure}

For any $\alpha, \beta \in \pi$, let $D_{-,-,+}(\alpha, \beta; 1, 1)$ be an unoriented $X$-cobordism depicted in Figure~$\ref{seki}$. 
The unoriented $X$-cobordism $D_{-,-,+}(\alpha, \beta; 1, 1)$ is a twice-punctured disk whose bottom base is the disjoint union of two unoriented $X$-manifolds $(\mathbf{S}^{1}, \alpha)$ and $(\mathbf{S}^{1}, \beta)$, whose top base is an unoriented $X$-manifold $(\mathbf{S}^{1}, \alpha\beta)$ and whose characteristic map sends the arcs labeled by $1\in\pi$ onto $x_{0}\in X$. 
For any  $\alpha, \beta \in \pi$, we define a homomorphism of $R$-modules $m_{\alpha, \beta}=m\colon A(\mathbf{S}^{1}, \alpha)\otimes A(\mathbf{S}^{1}, \beta)\rightarrow A(\mathbf{S}^{1}, \alpha\beta)$ by $\tau(D_{-,-,+}(\alpha, \beta; 1, 1))$. 
For any $v\in A(\mathbf{S}^{1}, \alpha)$ and $w\in A(\mathbf{S}^{1}, \beta )$, we denote $m(v\otimes w)$ by $vw\in A(\mathbf{S}^{1}, \alpha\beta)$. 

\begin{lem}\label{lem2}
For two elements $v\in A(\mathbf{S}^{1}, \alpha)$ and $w\in A(\mathbf{S}^{1}, \beta)$, we have $\Phi(vw)=\Phi(w)\Phi(v)$. 
\end{lem}
\begin{proof}
The proof of this lemma is similar to that of Lemma~$\ref{lem3}$. 
\end{proof}
\par
\begin{figure}[!h]
\begin{center}
\includegraphics[scale=0.26]{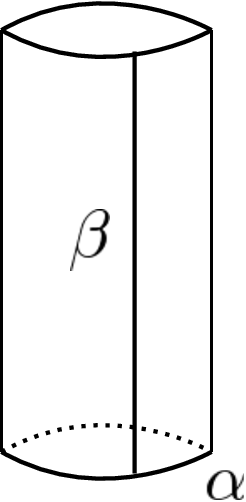}
\end{center}
\caption{The cobordism $C_{-,+}$($\alpha; \beta$).}
\label{phi}
\end{figure}
\begin{figure}[!h]
\begin{center}
\includegraphics[scale=0.425]{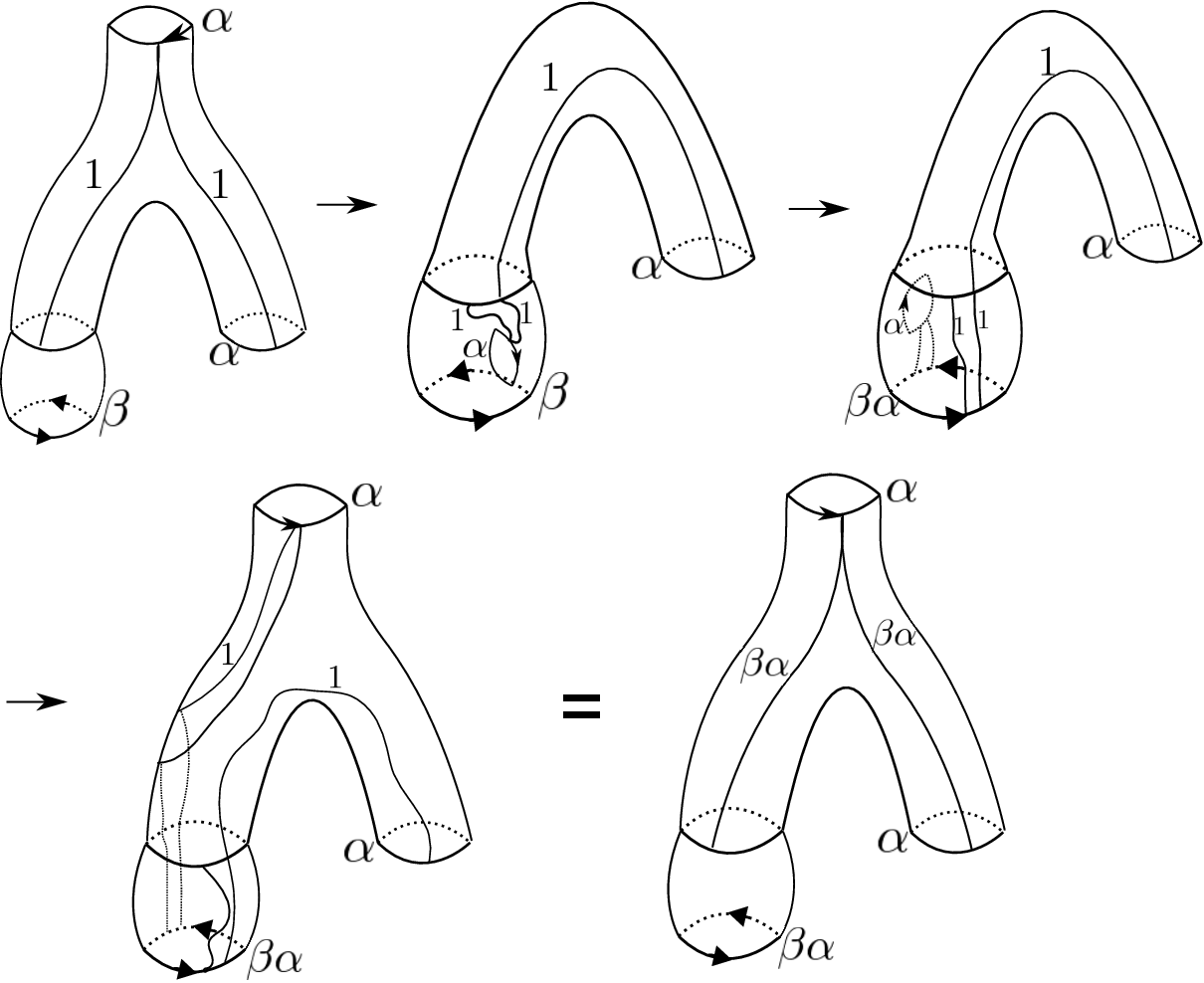}
\end{center}
\caption{Proof of  Lemma~$\ref{lem5}$.}
\label{proof5}
\end{figure}
\begin{figure}[!h]
\begin{center}
\includegraphics[scale=0.55]{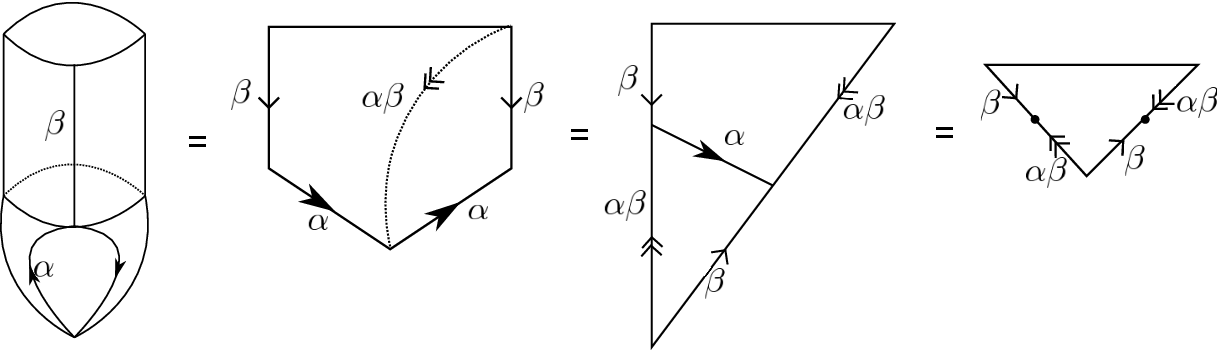}
\end{center}
\caption{Proof of Lemma~$\ref{lem6}$.}
\label{proof6}
\end{figure}
\begin{figure}[!h]
\begin{center}
\includegraphics[scale=0.6]{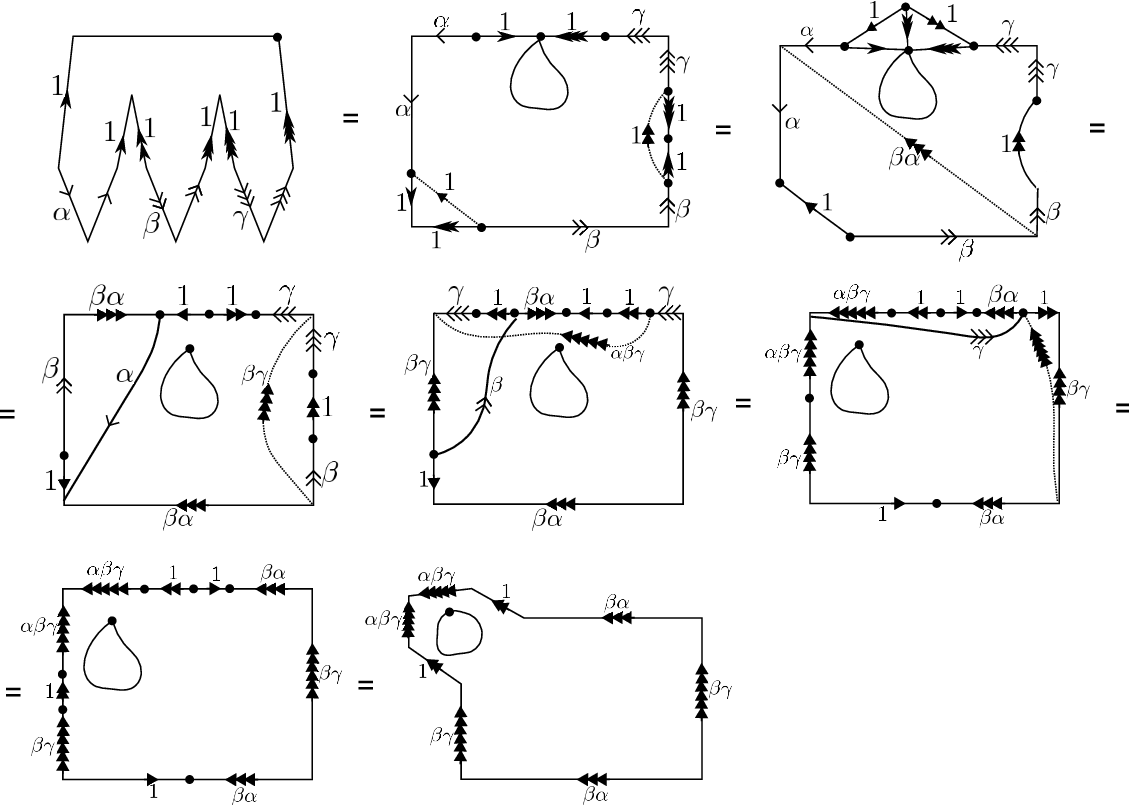}
\end{center}
\caption{Proof of Lemma~$\ref{lem7}$. As usual, we identify edges with some types of arrows in these pictures. }
\label{proof7}
\end{figure}
\begin{figure}[!h]
\begin{center}
\includegraphics[scale=0.4]{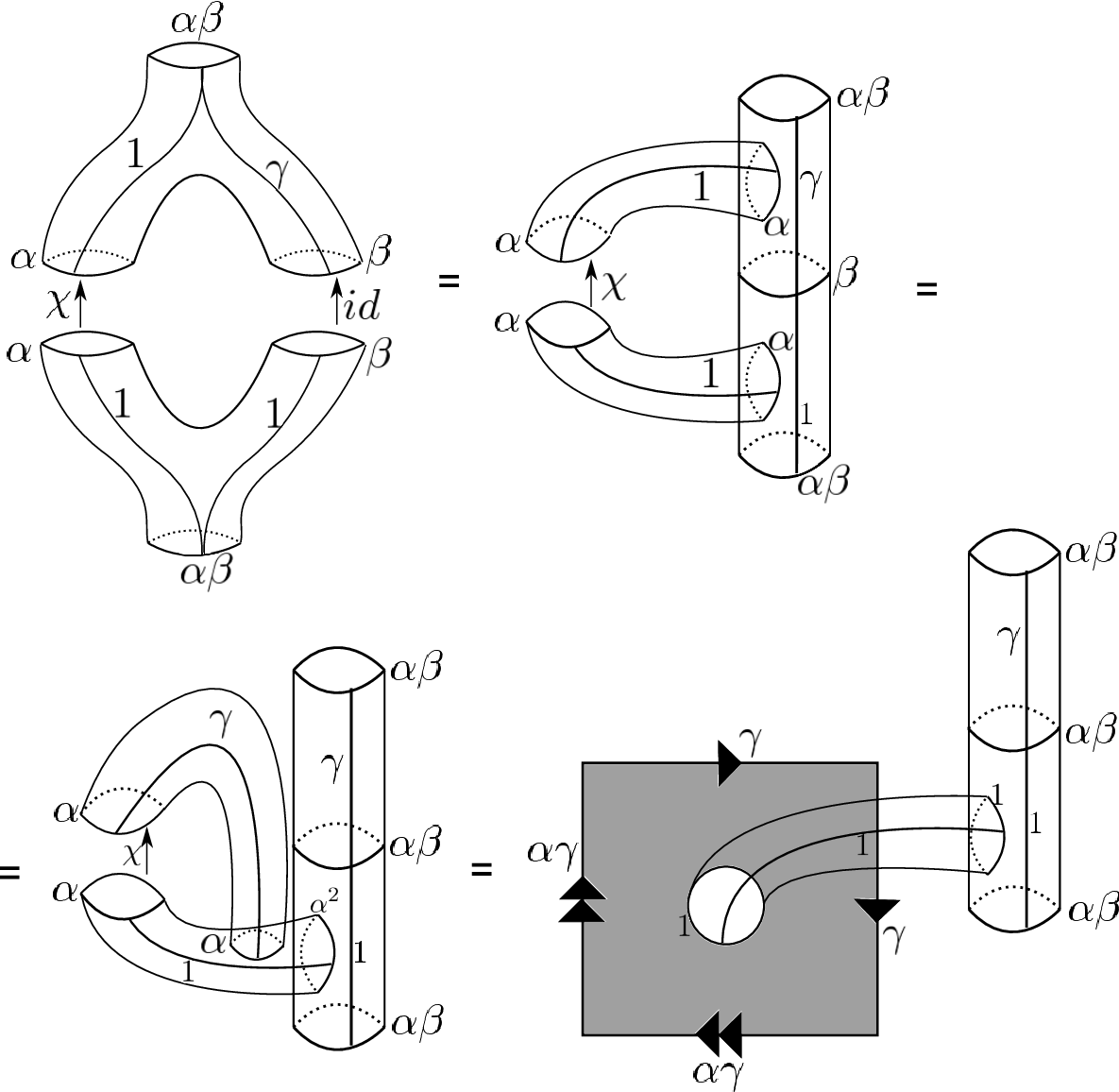}
\end{center}
\caption{Proof of Lemma~$\ref{lem8}$. $m\circ(\Phi \otimes \varphi _{\gamma })\circ\Delta _{\alpha, \beta}(v)=\varphi _{\gamma }(\theta_{\alpha\gamma }\theta_{\gamma }v)$. In the last picture, we identify edges with some types of arrows in these pictures.}
\label{tennkaizu3}
\end{figure}
\begin{figure}[!h]
\begin{center}
\includegraphics[scale=0.1]{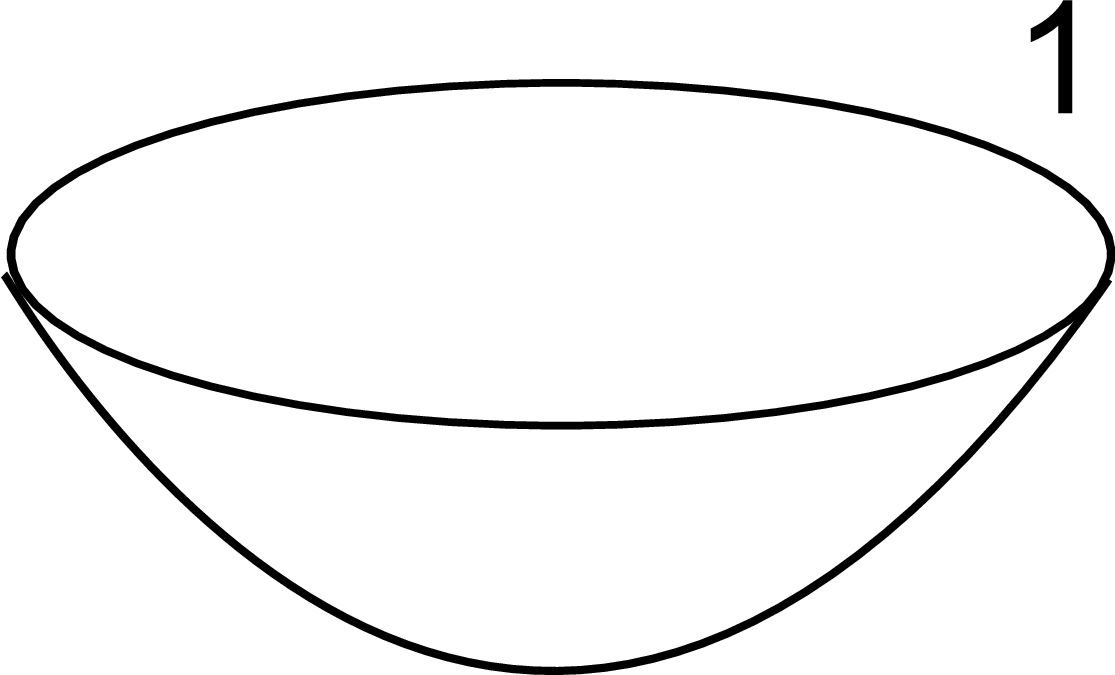}
\end{center}
\caption{The cobordism $(D, \emptyset, \partial{D})$.}
\label{1}
\end{figure}

For any $\alpha, \beta\in\pi$, let $C_{-,+}(\alpha; \beta)$ be an unoriented $X$-cobordism depicted in Figure~$\ref{phi}$. 
The unoriented $X$-cobordism $C_{-,+}(\alpha; \beta)$ is an annulus whose bottom base is an unoriented $X$-manifold $(\mathbf{S}^{1}, \alpha)$, whose top base is also an unoriented $X$-manifold $(\mathbf{S}^{1}, \alpha)$ and whose characteristic map sends the arc labeled by $\beta \in\pi$ onto a loop on $X$ whose homotopy class is $\beta \in\pi$. 
For any $\alpha, \beta\in \pi$, we define a homomorphism of $R$-modules $\varphi_{\beta}\colon \bigoplus_{\alpha \in \pi}A(\mathbf{S}^{1}, \alpha)\rightarrow \bigoplus_{\alpha \in \pi}A(\mathbf{S}^{1}, \alpha)$ by $\bigoplus_{\alpha \in \pi}\tau(C_{-,+}(\alpha; \beta))$. 
\begin{lem}\label{lem4}
For any $\beta\in\pi$, we have $\Phi\circ\varphi_{\beta}\circ\Phi=\varphi_{\beta}$. 
\end{lem}
\begin{proof}
We can prove this lemma by using the same argument as Lemma $\ref{lem3}$. 
\end{proof}
\begin{lem}\label{lem5}
For any $\alpha\in\pi$ and $v_{\alpha}\in A(\mathbf{S}^{1}, \alpha)$, we have $\Phi(\theta_{\beta}v_{\alpha})=\varphi_{\beta\alpha}(\theta_{\beta\alpha}v_{\alpha})$. 
\end{lem}
\begin{proof}
Figure $\ref{proof5}$ shows this lemma. 
In Figure $\ref{proof5}$ the first cobordism corresponds to $\theta_{\beta}v_{\alpha}$ and the fifth cobordism corresponds to $\varphi_{\beta\alpha}(\theta_{\beta\alpha}v_{\alpha})$, where two arrows depicted in Figure $\ref{proof5}$ mean that two edges endowed with the arrows are identified respecting the orientations indicated by them. 
Sliding the top base of the first cobordism, we obtain the second, the third and the fourth cobordisms. 
As a result the top base is reversed. 
From these transformations and Definition~$\ref{HQFT}$, we have $\Phi(\theta_{\beta}v_{\alpha})=\varphi_{\beta\alpha}(\theta_{\beta\alpha}v_{\alpha})$. 
\end{proof}
\begin{lem}\label{lem6}
For any $\alpha, \beta\in\pi$, we have $\varphi_{\beta}(\theta_{\alpha})=\theta_{\alpha}$. 
\end{lem}
\begin{proof}
In Figure~$\ref{proof6}$, the first cobordism corresponds to $\varphi_{\beta}(\theta_{\alpha})$, where arrows depicted in Figure~$\ref{proof6}$ mean that edges endowed with these arrows are identified along the same arrows. The fourth cobordism corresponds to $\theta_{\alpha}$ because $\beta\alpha\beta=\alpha$. 
\end{proof}
\begin{lem}\label{lem7}
For any $\alpha, \beta, \gamma \in\pi$, let 
$Q$ be the unoriented $(1+1)$-dimensional $X$-cobordism introduced in Definition~$\ref{extcross}$ and depicted in Figure~$\ref{tennkaizu1}$. 
Then we have $\theta_{\alpha }\theta_{\beta }\theta_{\gamma }=\tau(Q)(1)\theta_{\alpha \beta \gamma }$, where $1$ is the unit of $R$. 
\end{lem}
\begin{proof}
Figure~$\ref{proof7}$ shows this lemma. In Figure~$\ref{proof7}$ the first cobordism corresponds to $\theta_{\alpha }\theta_{\beta }\theta_{\gamma }$ and the eighth cobordism corresponds to $\tau(Q)(1)\theta_{\alpha \beta \gamma }$. 
\end{proof}
\begin{lem}\label{lem8}
For any $\alpha, \beta, \gamma \in\pi$ and $v\in A(\mathbf{S}^{1}, \alpha \beta )$, we have the following equations: 
\begin{center}
$m\circ(\Phi\otimes \varphi_{\gamma })\circ \Delta _{\alpha, \beta }(v)=\varphi_{\gamma }(\theta_{\alpha \gamma }\theta_{\gamma }v)$, 
\par
$m\circ(\varphi_{\gamma }\otimes \Phi)\circ \Delta _{\alpha, \beta }(v)=\varphi_{\gamma }(\theta_{\beta \gamma }\theta_{\gamma }v)$, 
\end{center}
where $\Delta _{\alpha, \beta }=\tau(D_{+,+,-}(\alpha, \beta; 1, 1))$ and 
$D_{+,+,-}(\alpha, \beta; 1, 1)$ is the unoriented $X$-cobordism introduced in Definition~$\ref{extcross}$.  
\end{lem}
\begin{proof}
This lemma follows from Figure~$\ref{tennkaizu3}$. In Figure~$\ref{tennkaizu3}$ the first cobordism corresponds to $m\circ(\varphi_{\gamma }\otimes \Phi)\circ \Delta _{\alpha, \beta }(v)$ and the fourth cobordism corresponds to $\varphi_{\gamma }(\theta_{\alpha \gamma }\theta_{\gamma }v)$. Similarly we can prove $m\circ(\varphi_{\gamma }\otimes \Phi)\circ \Delta _{\alpha, \beta }(v)=\varphi_{\gamma }(\theta_{\beta \gamma }\theta_{\gamma }v)$. 
\end{proof}
\begin{lem}\label{lem9}
We have $\Phi(1_{L})=1_{L}$, where $1_{L}=\tau(D, \emptyset, \partial{D})(1)$, $1$ is the unit of $R$ and $D$ is a cup which is an unoriented $X$-cobordism depicted in Figure~$\ref{1}$. Note that the characteristic map of $D$ is uniquely determined. 
\end{lem}
\begin{proof}
By using similar argument of Lemma~$\ref{lem3}$, we can prove this. 
\end{proof}

Theorem~$\ref{turaev}$ and Lemmas~$\ref{lem1}$-$\ref{lem9}$ show that any unoriented $(1+1)$-dimensional HQFT $(A, \tau)$ with target $X$ induces an extended crossed $\pi$-algebra. We call the extended crossed $\pi$-algebra the {\it underlying extended crossed $\pi$-algebra} of the unoriented $(1+1)$-dimensional HQFT $(A, \tau)$. 
\par
The following theorem is our main theorem which is proved in Sections~$\ref{main}$ and $\ref{main2}$. 
\begin{thm}[Main theorem]\label{mainthm}
Let $\pi$ be a group with $\alpha^{2}=1$ for any $\alpha\in\pi$ and $X$ be a $K(\pi, 1)$ space. Then every unoriented $(1+1)$-dimensional HQFT with target $X$ over the ring $R$ determines an underlying extended crossed $\pi$-algebra over $R$. This induces a bijection between the set of isomorphism classes of unoriented $(1+1)$-dimensional HQFTs over $R$ and the set of isomorphism classes of extended crossed $\pi$-algebras over $R$. 
\end{thm}
\section{Proof of main theorem}\label{main}
In this section, we prove Theorem~$\ref{mainthm}$. 
To prove the theorem, we need to make an unoriented HQFT $(A, \tau)$ from a given extended crossed $\pi$-algebra $(L$, $\eta$, $\varphi$, $\{\theta_{\alpha}\}_{\alpha\in\pi}$, $\Phi)$. 
Our proof has three steps. 
In Step $1$, we construct a functor $A$. 
In Step $2$, we make a functor $\tau$. 
In Step $3$, we prove that the pair $(A, \tau)$ satisfies the axioms of HQFTs. 
To construct them we use the same method as \cite{turner-turaev:2006}. 
\par
[Step $1$]: Construction of a functor $A$. 
\par
Let $(M, g)$ be a connected unoriented 1-dimensional $X$-manifold. We define an $R$-module $A(M, g)$ by 
\begin{align*}
A(M, g)
:=\{(r, v)|r\colon (\mathbf{S}^{1}, \alpha)\rightarrow (M, g): {\rm unoriented}\ X{\rm -homeomorphism}, v\in L_{\alpha}\}/\approx, 
\end{align*}
where $(r, v)\approx (r,' v')$ if and only if $r$ is isotopic to $r'$ and $v=v'$, or $r$ is not isotopic to $r'$ and $v=\Phi(v')$. For any $(M, g)$, such an $\alpha\in\pi$ is uniquely determined. For any unoriented  $X$-homeomorphism $h\colon (\mathbf{S}^{1}, \alpha)\rightarrow (M, g)$, we define a map $\widetilde{h}\colon A(M, g)\rightarrow L_{\alpha}$ by 
\begin{equation*}
\widetilde{h}(r, v):=
\begin{cases}
v   &\text{(if $r$ is isotopic to $h$)}, \\
\Phi (v)   &\text{(if $r$ is not isotopic to $h$)}. 
\end{cases}
\end{equation*}
Then the map $\widetilde{h}$ is bijective. In fact it has inverse map $(\widetilde{h})^{-1}\colon L_{\alpha}\rightarrow A(M, g)$ which is defined by $(\widetilde{h})^{-1}(v)=(h, v)$ for any $v\in L_{\alpha}$. 
Moreover we can use the $R$-module structure of $L_{\alpha}$ to turn $A(M, g)$ into an $R$-module. The $R$-module structure of $A(M,g)$ does not depend on the choice of the map $h$. 
This follows from the folowing: 
\begin{align*}
a(r,v)=& (\widetilde{h})^{-1}(a\widetilde{h}(r,v))\\
      =&
\begin{cases}
(h, av) &\text{(if $r$ is isotopic to $h$)}, \\
(h, a\Phi (v))&\text{(if $r$ is not isotopic to $h$)}\\
\end{cases}\\
      =&
\begin{cases}
(r, av) &\text{(if $r$ is isotopic to $h$)}, \\
(r, \Phi (a\Phi (v)))&\text{(if $r$ is not isotopic to $h$)}\\
\end{cases}\\
      =&(r, av), 
\end{align*} 
where $(r,v)\in A(M, g)$ and $a\in R$. 
Since $L_{\alpha}$ is a projective $R$-module, so is $A(M, g)$. 
In general we define $A(\emptyset)$ by $R$ and $A(M\sqcup N)$ by $A(M)\otimes A(N)$ for all connected unoriented 1-dimensional $X$-manifolds $M$ and $N$ (more precisely $M\sqcup N$ is an ordered disjoint union and $A(M)\otimes A(N)$ is an ordered tensor product). 
For any unoriented $X$-homeomorphism of unoriented $X$-manifolds $f\colon (M, g)\rightarrow (M', g')$, we define an $R$-homomorphism $f_{\sharp}\colon A(M, g)\rightarrow A(M', g')$ by $f_{\sharp}(r, v)=(f\circ r, v)$ for any $(r, v)\in\ A(M, g)$. 
\par
[Step $2$]: Construction of a functor $\tau$. 
\par
For any unoriented ($1+1$)-dimensional $X$-cobordism $(W, M_{0}, M_{1}, g)$, we define an $R$-homomorphism $\tau(W, g)\colon A(M_{0},g|_{M_{0}})\rightarrow A(M_{1},g|_{M_{1}})$ as follows: 
\par
Case $1$: $W$ is orientable and connected. 
\par
\begin{figure}[!h]
\begin{center}
\includegraphics[scale=0.35]{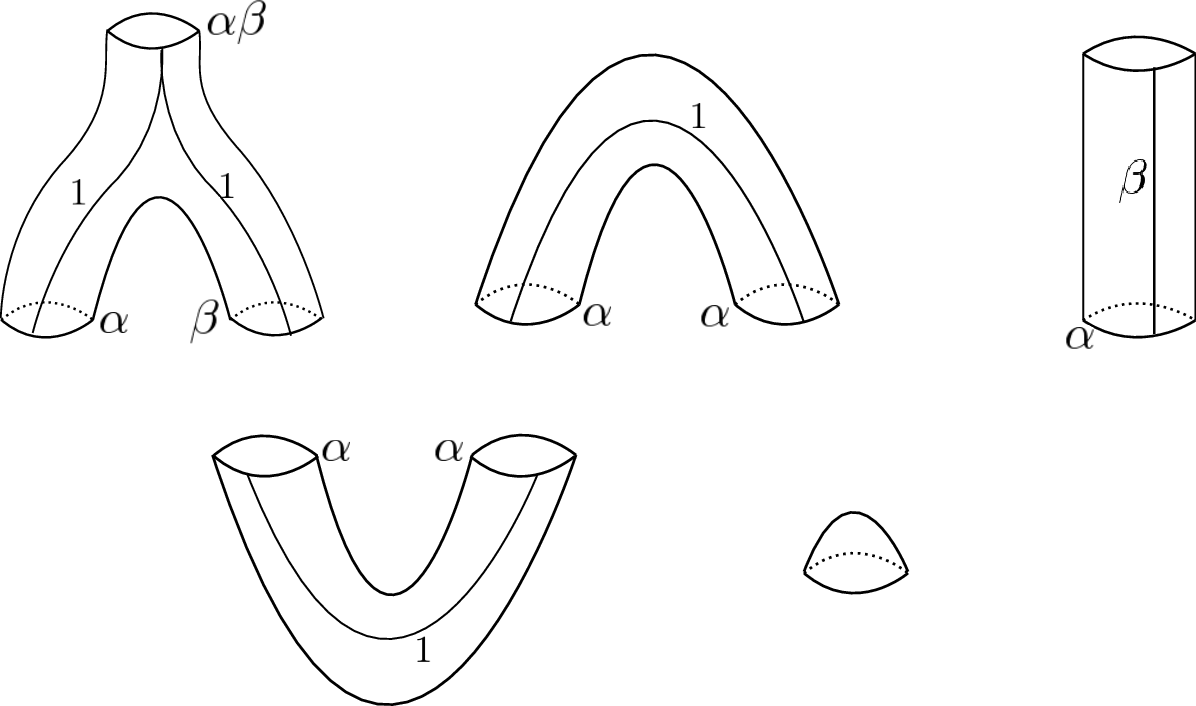}
\end{center}
\caption{Oriented basic cobordisms.}
\label{kihon cobordisms}
\end{figure}

Fix an orientation of $\mathbf{S}^{1}$ and give $W$ an orientation, and
take unoriented $X$-homeomorphisms $h_{M_{0}}\colon (\mathbf{S}^{1}, \alpha_{1})\sqcup\dots\sqcup(\mathbf{S}^{1}, \alpha_{n})\rightarrow (M_{0}, g|_{M_{0}})$ and $h_{M_{1}}\colon (\mathbf{S}^{1}, \beta_{1})\sqcup\dots\sqcup(\mathbf{S}^{1}, \beta_{n})\rightarrow (M_{1}, g|_{M_{1}})$ which preserve orientations. 
Then we define an $R$-homomorphism $\tau (W, g)\colon A(M_{0},g|_{M_{0}})\rightarrow A(M_{1},g|_{M_{1}})$ by $\widetilde{h_{M_{1}}}^{-1}\circ \tau^{L}(W,g)\circ \widetilde{h_{M_{0}}}$. 
The definition of $\tau^{L}$ is introduced in Theorem~$\ref{turaev}$. 
We need to prove that $\tau(W,g)$ does not depend on the choice of their orientations. 
It is sufficient that we check it in the cases where $W$ is an unoriented basic cobordism depicted in Figure~$\ref{kihon cobordisms}$. 
When $W$ is an unorieted $X$-cobordism at the upper left in Figure~$\ref{kihon cobordisms}$, 
take unoriented $X$-homeomorphisms $h_{M_{0}}\colon (\mathbf{S}^{1}, \alpha)\sqcup(\mathbf{S}^{1}, \beta )\rightarrow M_{0}$ and $h_{M_{1}}\colon (\mathbf{S}^{1}, \alpha\beta )\rightarrow M_{1}$. Then we have 
\begin{align*}
\widetilde{h_{M_{1}}}^{-1}\circ m\circ\widetilde{h_{M_{0}}}&= \widetilde{h_{M_{1}}}^{-1}\circ \Phi \circ m\circ P \circ (\Phi \otimes\Phi )\circ\widetilde{h_{M_{0}}}\\
&=(\widetilde{h_{M_{1}}\circ \chi })^{-1}\circ m\circ P \circ (\widetilde{h_{M_{0}}\circ (\chi \sqcup\chi ))}, 
\end{align*}
where $P$ is the permutation. 
This equation implies that $\tau(W, g)$ does not depend on the choice of the orientation of $W$. 
In other cases, we can use similar arguments. 
\par
Case $2$: $W$ is non-orientable and connected.  
\par
\begin{figure}[!h]
\begin{center}
\includegraphics[scale=0.23]{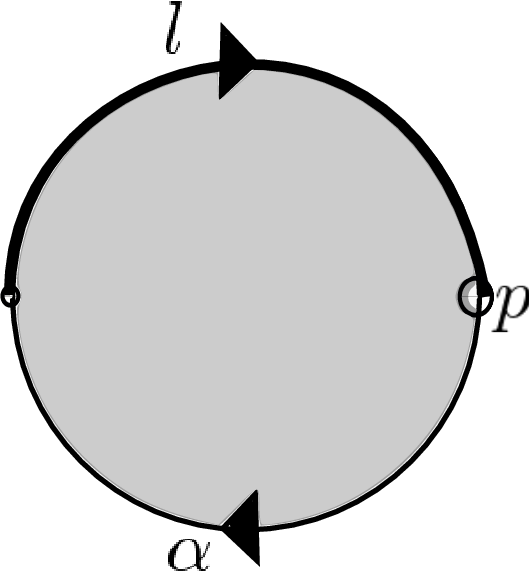}
\end{center}
\caption{The cobordism $(\mathbf{R}P^{2}, f_{\alpha })$D}
\label{tennkaizu4}
\end{figure}

Let $\mathbf{R}P^{2}$ be the projective plane. 
For any $\alpha\in\pi$, we define an unoriented $X$-cobordism $(\mathbf{R}P^{2}, f_{\alpha}, p)$ with $p$ a point of $\mathbf{R}P^{2}$ as follows. 
The pair $(\mathbf{R}P^{2}, f_{\alpha})$ is an unorieted $X$-cobordism $(\mathbf{R}P^{2}, \emptyset, \emptyset, f_{\alpha})$ such that $f_{\alpha}(p)=x_{0}$ and that the homotopy class of $f_{\alpha}|_{l}$ equals $\alpha\in\pi$ for the loop $l$ on $\mathbf{R}P^{2}$ depicted in Figure~$\ref{tennkaizu4}$
(in Figure~$\ref{tennkaizu4}$, $l$ is presented by the upper arc with arrow, which is identified with the lower arc with arrow).
Such an unoriented $X$-cobordism is uniquely determined up to homotopy by $\alpha:=[l]\in\pi$.  
In general for any unoriented $X$-cobordism $(\mathbf{R}P^{2}, g)$, (by using a homotopy) 
we can assume that there are $p\in\mathbf{R}P^{2}$ and $\alpha\in\pi$ which satisfy the following conditions:
\begin{itemize}
\item $g(p)=x_{0}$, 
\item $g$ is homotopic to $f_{\alpha}$.  
\end{itemize}
Then we denote the unoriented $X$-cobordism $(\mathbf{R}P^{2}, g)$ by $(\mathbf{R}P^{2}, f_{\alpha})$ . 
Now we represent $(W, g)$ as the connected sum of an orientable (unoriented) $X$-cobordism $(W^{{\rm or}}, g|_{W^{{\rm or}}})$ and unoriented $X$-cobordisms $(\mathbf{R}P^{2}, f_{\alpha _{1}}), \dots , (\mathbf{R}P^{2}, f_{\alpha _{n}})$, that is, 
\begin{center}
$(W, g)\cong (W^{{\rm or}}, g|_{W^{{\rm or}}})\sharp(\mathbf{R}P^{2}, f_{\alpha _{1}})\sharp\dots \sharp(\mathbf{R}P^{2}, f_{\alpha _{n}})$. 
\end{center}
Note that $\partial{W^{{\rm or}}}=\partial{W}$ and that a homomorphism $\tau(W^{{\rm or}}, g|_{W^{{\rm or}}})$ is defined in the orientable case. 
Let $m$ be the number of components of $M_{1}$. We define a homomorphism $\tau(W, g)$ as follows. 
If $m>0$, take an unoriented $X$-homeomorphism $h\colon (\mathbf{S}^{1}, \beta_{1})\sqcup\cdots\sqcup(\mathbf{S}^{1}, \beta_{m})\rightarrow (M_{1}, g|_{M_{1}})$ and identify two $R$-modules $A(M_{1}, g|_{M_{1}})$ and $\bigotimes_{i=1}^{m}L_{\beta _{i}}$ by $\tilde{h}\colon A(M_{1}, g|_{M_{1}})\rightarrow \bigotimes_{i=1}^{m}L_{\beta _{i}}$. 
Under this identification, we define a map $\psi_{\alpha _{1}, \cdots, \alpha _{n}}\colon A(M_{1}, g|_{M_{1}})\rightarrow A(M_{1}, g|_{M_{1}})$ to be the identity on all factors except one where it is multiplication by $\prod_{i=1}^{n}\theta _{\alpha_{i}}$. We define 
\begin{center}
$\tau(W, g):=\psi_{\alpha _{1}, \cdots, \alpha _{n}}\circ \tau(W^{{\rm or}}, g|_{W^{{\rm or}}})$. 
\end{center}
If $m=0$, consider an unoriented $X$-cobordism $W^{{\rm or}}-D^{2}=(W^{{\rm or}}-D^{2}, M_{0}, M_{1}\sqcup \partial(D^{2}), g|_{W^{{\rm or}}-D^{2}})$, where $D^{2}$ is any disk on $W$. 
Take an unoriented $X$-homeomorphism $h\colon (\mathbf{S}^{1}, 1)\rightarrow (\partial(D^{2}), g|_{\partial(D^{2})})$ and identify two $R$-modules $A(\partial(D^{2}), g|_{\partial(D^{2})})$ and $L_{1}$ by $\tilde{h}\colon A(\partial(D^{2}), g|_{\partial(D^{2})})\rightarrow L_{1}$. 
Under this identification, we define a homomorphism $\tau(W, g)$ by 
\begin{center}
$\tau(W, g)(v):=\eta (\tau(W^{{\rm or}}-D^{2}, g|_{W^{{\rm or}}-D^{2}})(v),  \displaystyle{\prod_{i=1}^{n}\theta _{\alpha_{i}}})$
\end{center}
for any $v\in A(M_{0}, g|_{M_{0}})$. 
From Lemmas~$\ref{claim1}$ and $\ref{claim}$ below, the functor $\tau$ is well defined. 
\par
Case $3$: $W$ is not connected:
\par
We can extend the definition of $\tau$ constructed as above to non-connected cases by using tensor products as Step$1$. 
\par
[Step $3$]: The pair $(A, \tau)$ is an unoriented ($1+1$)-dimensional HQFT with target $X$. 
\par
We need to check the axioms of unoriented HQFTs (see Definition $\ref{HQFT}$). The pair $(A, \tau)$ clearly satisfies the axioms except for ($4$) and ($6$). 
In the next section, we show that  $(A, \tau)$ satisfies the axioms ($4$) and ($6$) (Propositions~$\ref{natural}$ and $\ref{gluing}$). 
\par
\begin{figure}[!h]
\begin{center}
\includegraphics[scale=0.45]{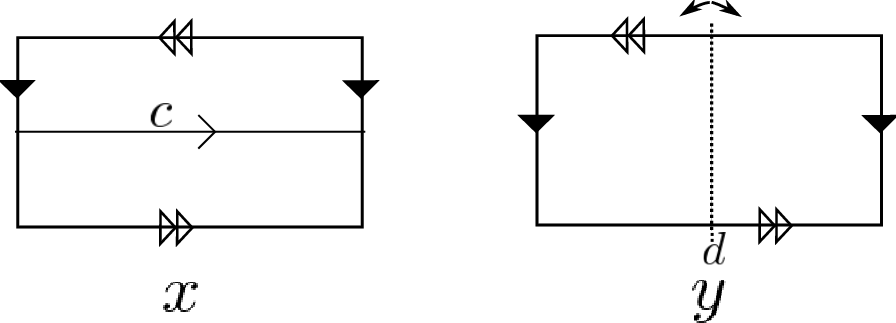}
\end{center}
\caption{Two generators of the mapping class group of the Klein bottle.}
\label{klein}
\end{figure}

From Step $1$, $2$ and $3$, we complete the proof of Theorem $\ref{mainthm}$ (except for Propositions~$\ref{natural}$ and $\ref{gluing}$ and Lemmas~$\ref{claim1}$ and $\ref{claim}$). 
\begin{lem}\label{claim1}
{\rm(i)} The map $\tau(W, g)$ does not depend on the choice of $h$. 
\par
{\rm(ii)} The map $\tau(W, g)$ does not depend on the choice of a factor multiplied the element $ \prod_{i=1}^{n}\theta _{\alpha_{i}} \in L_{1}$. 
\par
{\rm(iii)} The map $\tau(W, g)$ does not depend on the choice of the connected sum $(W, g)\cong (W^{{\rm or}}, g|_{W^{{\rm or}}})\sharp(\mathbf{R}P^{2}, f_{\alpha _{1}})\sharp\dots \sharp(\mathbf{R}P^{2}, f_{\alpha _{n}})$. 
\end{lem}
\begin{proof}
(i): In the case where $m=1$, take any unoriented $X$-homeomorphism $h\colon (\mathbf{S}^{1}, \alpha)\rightarrow (M_{1}, g|_{M_{1}})$. For any $(r, v)\in\ A(M_{1}, g|_{M_{1}})$, we have 
\begin{align*}
(\tilde{h})^{-1}(\displaystyle{\prod_{i=1}^{n}\theta _{\alpha_{i}}}\tilde{h}(r,v))
=&
\begin{cases}
(r, \displaystyle{\prod_{i=1}^{n}\theta _{\alpha_{i}}}v) &\text{(if $r$ is isotopic to $h$)},\\
(r, \Phi (\displaystyle{\prod_{i=1}^{n}\theta _{\alpha_{i}}}\Phi (v))) &\text{(if $r$ is not isotopic to $h$)}\\
\end{cases}\\
      =&(r, \displaystyle{\prod_{i=1}^{n}\theta _{\alpha_{i}}}v). 
\end{align*}
Hence the map $\tau(W, g)$ does not depend on the choice of $h$. Similarly we can prove the case where $m>1$. In the case where $m=0$, (i) follows from the fact that $\Phi$ preserves $\eta$ and $\theta_{\alpha}$ for any $\alpha$. 
\par
(ii): It follows from Figure~$\ref{kakeru2}$ and Theorem~$\ref{turaev}$. 
\begin{figure}[!h]
\begin{center}
\includegraphics[scale=0.45]{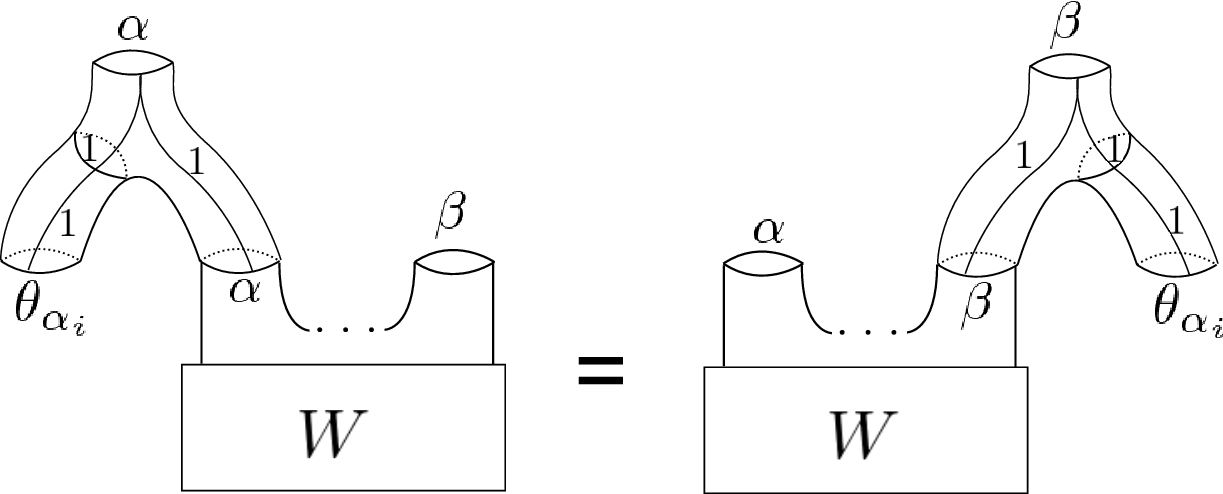}
\end{center}
\caption{Cobordism relation. }
\label{kakeru2}
\end{figure}
\par
(iii): In the proof of Lemma~$\ref{lem7}$, we proved that $(\mathbf{R}P^{2}, f_{\alpha})\sharp(\mathbf{R}P^{2}, f_{\beta})\sharp(\mathbf{R}P^{2}, f_{\gamma} )$ is unoriented $X$-homeomorphic to $(T^{2}, g_{\alpha , \beta , \gamma })\sharp(\mathbf{R}P^{2}, f_{\alpha\beta \gamma} )$, where $(T^{2}, g_{\alpha , \beta , \gamma })$ is the unoriented $X$-cobordism depicted in Figure~$\ref{3multi}$ whose bottom base and top base are empty and whose characteristic map $g_{\alpha , \beta , \gamma }$ sends the arcs labeled by $1, \beta \alpha , \beta \gamma \in\pi$ onto the loops with the corresponding labels. 
It follows from the definitions that $\tau((\mathbf{R}P^{2}, f_{\alpha})\sharp(\mathbf{R}P^{2}, f_{\beta})\sharp(\mathbf{R}P^{2}, f_{\gamma} ))=\tau((T^{2}, g_{\alpha , \beta , \gamma })\sharp(\mathbf{R}P^{2}, f_{\alpha\beta \gamma }))$. 
Hence it is sufficient to prove (iii) for the case where $n=1$ or $2$, which is shown in Lemma~$\ref{claim}$ below. 
\begin{figure}[!h]
\begin{center}
\includegraphics[scale=0.3]{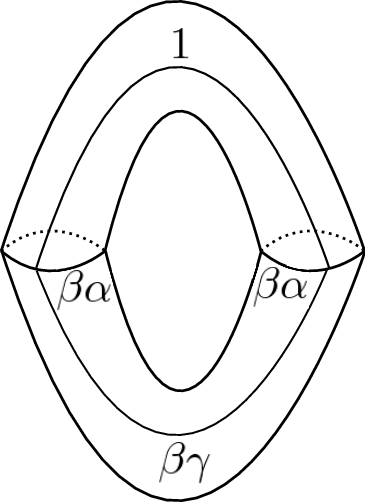}
\end{center}
\caption{The cobordism $(T^{2}, g_{\alpha , \beta , \gamma })$D}
\label{3multi}
\end{figure}
\end{proof}
\begin{lem}\label{claim}
{\rm(I)} Assume that we have two connected sums 
\begin{center}
$(W, g)\cong (W^{{\rm or}}, g|_{W^{{\rm or}}})\sharp(\mathbf{R}P^{2}, f_{\alpha })$ 
\end{center}
and 
\begin{center}
$(W, g)\cong(\widetilde{W}^{{\rm or}}, g|_{\widetilde{W}^{{\rm or}}})\sharp(\mathbf{R}P^{2}, f_{\beta} )$. 
\end{center}
\par
{\rm(I-a)} If we have an unoriented $X$-homeomorphism of unoriented $X$-manifolds $f\colon (W^{{\rm or}}, g|_{W^{{\rm or}}})\rightarrow (\widetilde{W}^{{\rm or}}, g|_{\widetilde{W}^{{\rm or}}})$ and $\alpha =\beta $, we have $\psi_{\alpha }\circ \tau(W^{{\rm or}}, g|_{W^{{\rm or}}})=\psi_{\beta}\circ \tau(\widetilde{W}^{{\rm or}}, g|_{\widetilde{W}^{{\rm or}}})$. 
\par
{\rm(I-b)} If we have an unoriented $X$-homeomorphism of unoriented $X$-manifolds $f\colon (\mathbf{R}P^{2}, f_{\alpha})\rightarrow (\mathbf{R}P^{2}, f_{\beta})$ and $W^{{\rm or}}=\widetilde{W}^{{\rm or}}$, we have $\psi_{\alpha}\circ \tau(W^{{\rm or}}, g|_{W^{{\rm or}}})=\psi_{\beta}\circ \tau(\widetilde{W}^{{\rm or}}, g|_{\widetilde{W}^{{\rm or}}})$. 
\par
{\rm(I-c)} We have $\psi _{ \alpha }\circ \tau(W^{{\rm or}}, g|_{W^{{\rm or}}})=\psi_{\beta}\circ \tau(\widetilde{W}^{{\rm or}}, g|_{\widetilde{W}^{{\rm or}}})$. 
\par
{\rm(II)} Assume that we have two connected sums 
\begin{center}
$(W, g)\cong (W^{{\rm or}}, g|_{W^{{\rm or}}})\sharp(\mathbf{R}P^{2}, f_{\alpha_{1}})\sharp(\mathbf{R}P^{2}, f_{\alpha_{2}})$ 
\end{center}
and 
\begin{center}
$(W, g)\cong(\widetilde{W}^{{\rm or}}, g|_{\widetilde{W}^{{\rm or}}})\sharp(\mathbf{R}P^{2}, f_{\beta_{1}} )\sharp(\mathbf{R}P^{2}, f_{\beta_{2}})$. 
\end{center}
\par
{\rm(II-a)} If we have an unoriented $X$-homeomorphism of unoriented $X$-manifolds $f\colon (W^{{\rm or}}, g|_{W^{{\rm or}}})\rightarrow (\widetilde{W}^{{\rm or}}, g|_{\widetilde{W}^{{\rm or}}})$ and $\{\alpha_{1}, \alpha_{2}\} =\{\beta_{1}, \beta_{2}\} $, we have $\psi_{ \alpha_{1}, \alpha _{2} }\circ \tau(W^{{\rm or}}, g|_{W^{{\rm or}}})=\psi_{\beta_{1}, \beta _{2}}\circ \tau(\widetilde{W}^{{\rm or}}, g|_{\widetilde{W}^{{\rm or}}})$. 
\par
{\rm(II-b)} If we have an unoriented $X$-homeomorphism of unoriented $X$-manifolds $f\colon (\mathbf{R}P^{2}, f_{\alpha_{1}})\sharp(\mathbf{R}P^{2}, f_{\alpha_{2}})\rightarrow (\mathbf{R}P^{2}, f_{\beta_{1}})\sharp(\mathbf{R}P^{2}, f_{\beta_{2}})$ and $W^{{\rm or}}=\widetilde{W}^{{\rm or}}$, we have $\psi_{\alpha_{1}, \alpha _{2}}\circ \tau(W^{{\rm or}}, g|_{W^{{\rm or}}})=\psi_{\beta_{1}, \beta _{2}}\circ \tau(\widetilde{W}^{{\rm or}}, g|_{\widetilde{W}^{{\rm or}}})$. 
\par
{\rm(II-c)}We have $\psi_{\alpha_{1}, \alpha _{2}}\circ \tau(W^{{\rm or}}, g|_{W^{{\rm or}}})=\psi_{\beta_{1}, \beta _{2}}\circ \tau(\widetilde{W}^{{\rm or}}, g|_{\widetilde{W}^{{\rm or}}})$. 
\end{lem}
\begin{proof}
(I-a): We can naturally identify $\partial(W^{{\rm or}})$ and $\partial(\widetilde{W}^{{\rm or}})$ with $\partial(W)$. It follows from the definition of $\tau(W^{{\rm or}})$ and $\tau(\widetilde{W}^{{\rm or}})$ and Theorem~$\ref{HQFT}$ that 
\begin{align}
(f|_{M_{1}})_{\sharp}\circ\tau(W^{{\rm or}})=\tau(\widetilde{W}^{{\rm or}})\circ(f|_{M_{0}})_{\sharp}. \label{a1}
\end{align}
The map $(f|_{M_{1}})_{\sharp}$ is the identity map or $\Phi$ on each factor of $A(M_{1}, g|_{M_{1}})$. 
It follows from the definition of extended crossed $\pi$-algebras that the center of $L$ contains $L_{1} $ and that $\Phi (\theta _{\alpha }v)=\theta _{\alpha }\Phi (v)$ for any $\alpha\in\pi$ and $v\in L$. 
Hence we have 
\begin{align}
\psi_{\alpha}\circ(f|_{M_{1}})_{\sharp}=(f|_{M_{1}})_{\sharp}\circ\psi_{\alpha}. \label{a2}
\end{align}
It follows from ($\ref{a1}$) and ($\ref{a2}$) that 
\begin{align}
(f|_{M_{1}})_{\sharp}\circ\psi_{\alpha}\circ\tau(W^{{\rm or}})=\psi_{\alpha}\circ\tau(\widetilde{W}^{{\rm or}})\circ(f|_{M_{0}})_{\sharp}. \label{a3}
\end{align}
The equation ($\ref{a3}$) implies (I-a). 
\par
(I-b): 
Since the mapping class group of the projective plane $\mathbf{R}P^{2}$ is trivial, the map $f$ is isotopic to the identity map, that is, there exists a continuous map $H\colon \mathbf{R}P^{2}\times [0,1]\rightarrow \mathbf{R}P^{2}$ such that the map $H_{t}\colon \mathbf{R}P^{2}\rightarrow \mathbf{R}P^{2}$ defined by $H_{t}(x)=H(x, t)$ for any $x\in\mathbf{R}P^{2}$ is a homeomorphism for all $t\in [0, 1]$ and satisfies $H_{0}=f$ and $H_{1}=\id$. 
We do not know if the map $H$ fixes $p\in\mathbf{R}P^{2}$, where $p$ is the point introduced in the proof of Theorem~$\ref{mainthm}$. 
Let $\gamma \in\pi$ be an element corresponding to the loop $f_{\beta}(H(\{p\}\times [0, 1]))$ on $X$. Then we have $\alpha=\gamma \beta \gamma =\beta \gamma ^{2}=\beta $. This implies (I-b). 
\par
(I-c): Consider the loops $c$, $f(c)$ and $c'$ as depicted in Figure~$\ref{proof-c}$. 
In general, $f(c)$ is not homotopic to $c'$. 
For any $\alpha, \beta\in\pi$ and $v\in L_{\alpha}$, we have $\Phi(\theta_{\beta}v)=\varphi_{\alpha \beta }(\theta_{\alpha \beta }v)$. 
This equation means that $\tau(W, g)$ is preserved under the transformation depicted in Figure~$\ref{henkei}$ (see the proof of Lemma~$\ref{lem5}$). 
By the definition of $\tau(W, g)$ and Theorem~$\ref{turaev}$, the map $\tau(W, g)$ is preserved by Dehn twists on $W$. 
By using the transformation depicted in Figure~$\ref{henkei}$ and Dehn twists, we can assume that $f(c)$ is homotopic to $c'$. 
If $f(c)$ is homotopic to $c'$, there exists an unoriented $X$-homeomorphism $f'\colon (\mathbf{R}P^{2}, f_{\alpha})\rightarrow (\mathbf{R}P^{2}, f_{\beta})$ and we can use arguments in (I-b). 
\begin{figure}[!h]
\begin{center}
\includegraphics[scale=0.5]{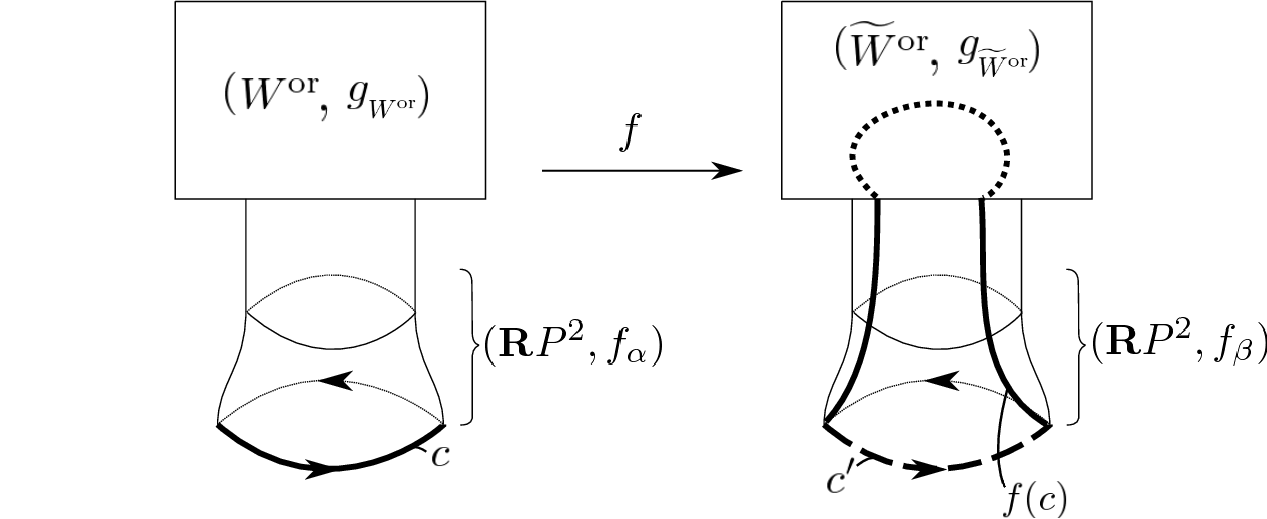}
\end{center}
\caption{Loops $c, f(c)$ and $c'$. }
\label{proof-c}
\end{figure}
\begin{figure}[!h]
\begin{center}
\includegraphics[scale=0.5]{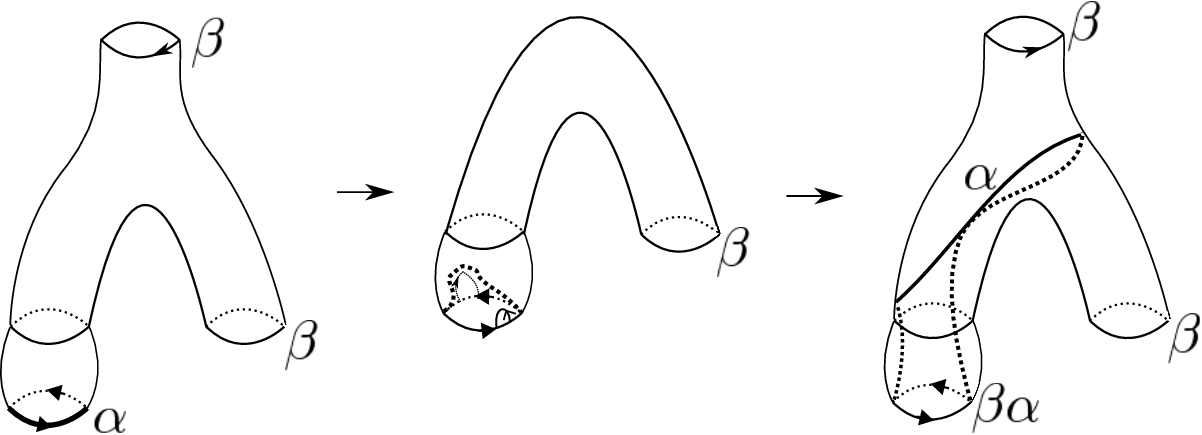}
\end{center}
\caption{A transformation. }
\label{henkei}
\end{figure}
\par
(II-a): 
We can show this by using similar arguments in the proof of (I-a). 
\par
(II-b): 
We can show this by using similar arguments in the proof of (I-b). Instead of the mapping class group of $\mathbf{R}P^{2}$, we use that of the Klein bottle. The mapping class group of a Klein bottle is generated by two elements $x$ and $y$ (see Theorem~$\ref{mcg-of-klein}$ below). 
We can assume that the cobordism $(\mathbf{R}P^{2}, f_{\alpha} )\sharp(\mathbf{R}P^{2}, f_{\beta} )$ is given by the right hand side in Figure~$\ref{tennkaizu5}$. 
If $f$ is isotopic to $x$, we have $\alpha _{1}=\beta _{2}$ and $\alpha _{2}=\beta _{1}$. If $f$ is isotopic to $y$, we have $\alpha _{1}=\beta _{1}$, $\alpha _{2}=\beta _{2}$ (see Figures~$\ref{x}$ and $\ref{y}$).
\par
(II-c): We can show this by using similar arguments in the proof of (I-c). 
\begin{figure}[!h]
\begin{center}
\includegraphics[scale=0.5]{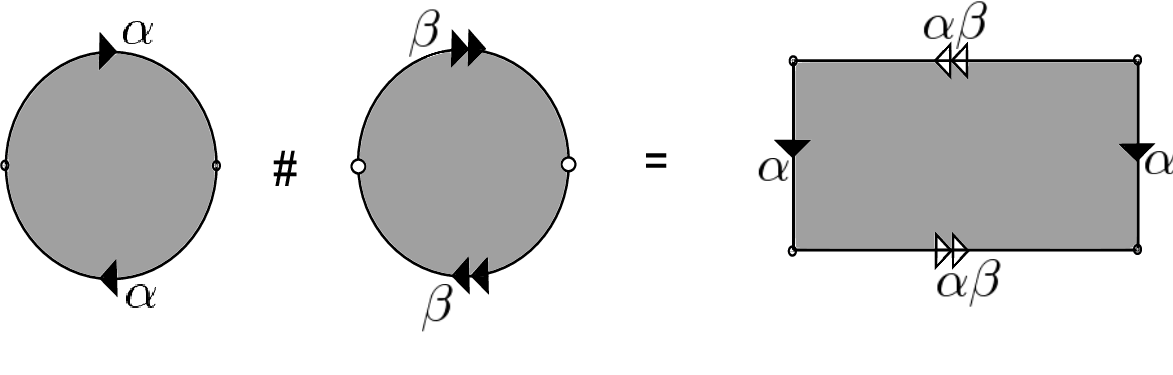}
\end{center}
\caption{The cobordism $(\mathbf{R}P^{2}, f_{\alpha} )\sharp(\mathbf{R}P^{2}, f_{\beta} )$.}
\label{tennkaizu5}
\end{figure}
\begin{figure}[!h]
\begin{center}
\includegraphics[scale=0.4]{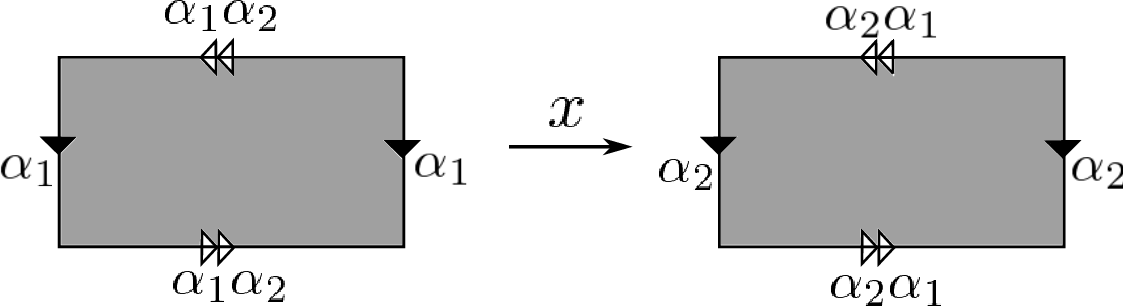}
\end{center}
\caption{In the case $f=x$. }
\label{x}
\end{figure}
\begin{figure}[!h]
\begin{center}
\includegraphics[scale=0.4]{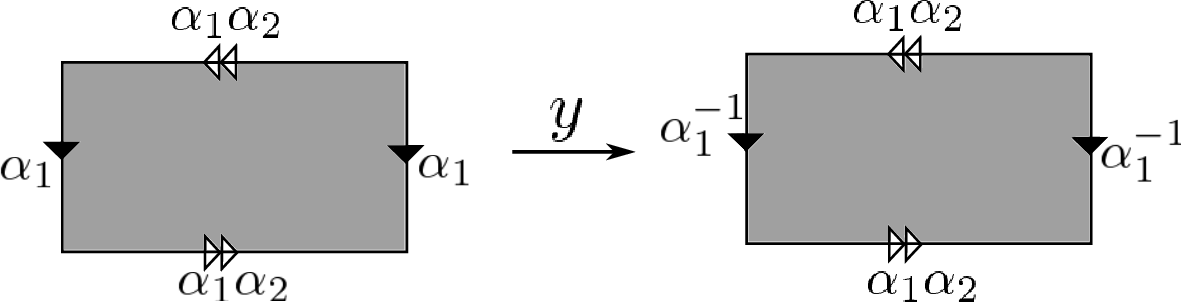}
\end{center}
\caption{In the case $f=y$. }
\label{y}
\end{figure}
\end{proof}
\begin{thm}[\cite{mcg-of-klein}]\label{mcg-of-klein}
Let K be the Klein bottle. We define a homeomorphism $x\colon K\rightarrow K$ as a Dehn twist along the loop $c$ depicted in Figure~$\ref{klein}$ and a homeomorphism $y\colon K\rightarrow K$ as taking the mirror image with respect to the line $d$ depicted in Figure~$\ref{klein}$. Then the mapping class group of $K$ is generated by the isotopy classes of $x$ and $y$. 
\end{thm}

\section{The axioms of HQFT}\label{main2}
\begin{figure}[!h]
\begin{center}
\includegraphics[scale=0.5]{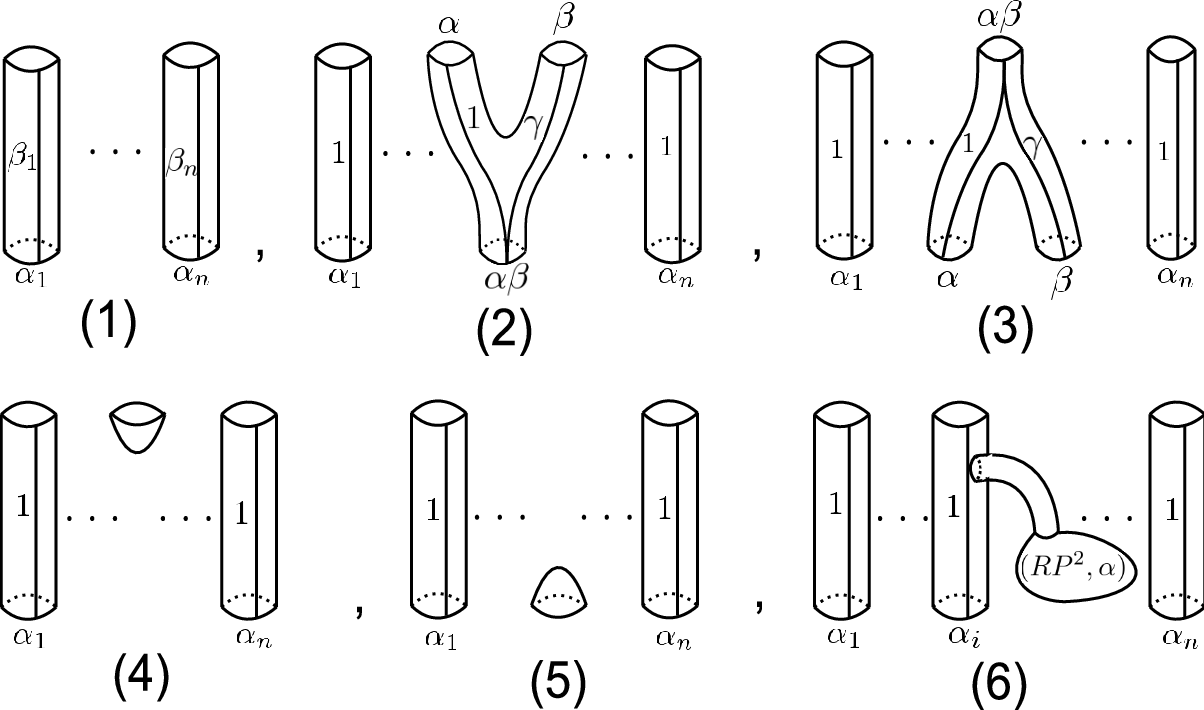}
\end{center}
\caption{Basic X-cobordisms.}
\label{basic cobordism}
\end{figure}

In this section, we check that the pair ($A, \tau$) constructed in Section~$\ref{main}$ satisfies the axioms of HQFTs (see Definition~$\ref{HQFT}$). 
\begin{prop}\label{natural}
The functor $\tau$ constructed in Theorem~$\ref{mainthm}$ from an extended crossed $\pi$-algebra $(L, \eta, \varphi, \{\theta_{\alpha}\}_{\alpha\in\pi}, \Phi)$ is natural with respect to unoriented $X$-homeomorphisms of unoriented $X$-manifolds. 
\end{prop}
\begin{proof}
Let $(W_{1}, M_{1}, N_{1}, g_{1})$ and $(W_{2}, M_{2}, N_{2}, g_{2})$ be two unoriented X-cobordisms and $f\colon (W_{1}, M_{1}, N_{1}, g_{1})\rightarrow (W_{2}, M_{2}, N_{2}, g_{2})$ be an unoriented X-homeomorphism of unoriented $X$-cobordisms. 
If we are given an unoriented $X$-homeomorphism  
\begin{center}
$(W_{1}, g_{1})\cong (W_{1}^{or}, g_{1}|_{W_{1}^{{\rm or}}})\sharp(\mathbf{R}P^{2}, f_{\alpha_{1} })\sharp\dots\sharp(\mathbf{R}P^{2}, f_{\alpha_{n} })$, 
\end{center}
then we have
\begin{center}
$(W_{2}, g_{2})\cong (f(W_{1}^{{\rm or}}), (g_{1}|_{W_{1}^{{\rm or}}})\circ 
f^{-1})\sharp(\mathbf{R}P^{2}, f_{\alpha_{1} }\circ f^{-1})\sharp\dots\sharp(\mathbf{R}P^{2}, f_{\alpha_{n} }\circ f^{-1})$. 
\end{center}
There is an element $\beta_{i}\in\pi$ such that $f_{\beta_{i}}$ is isotopic to $f_{\alpha_{i}}\circ f^{-1}$ for all $i=1,\dots, n$. 
It follows from the proof of Lemma~$\ref{claim}$ (I-b) that $\alpha_{i}=\beta_{i}$ for all $i=1,\dots, n$. From Theorem~$\ref{turaev}$, we have $(f|_{_{N_{1}}})_{\sharp}\circ\tau(W_{1}^{{\rm or}}, g_{1}|_{W_{1}^{{\rm or}}})=\tau(f(W_{1}^{{\rm or}}), (g_{1}|_{W_{1}^{{\rm or}}})\circ 
f^{-1})\circ(f|_{_{M_{1}}})_{\sharp}$. This completes the proof. 
\end{proof}
\begin{defn}
Let $(A, \tau)$ be the pair constructed in Section~$\ref{main}$ from an extended crossed $\pi$-algebra $(L, \eta, \varphi, \{\theta_{\alpha}\}_{\alpha\in\pi}, \Phi)$. 
An unoriented $X$-cobordism $(W_{0}$, $M_{0}$, $N)$ is {\it $X$-nice} if for any unoriented $X$-cobordism $(W_{1}, N' ,M_{1})$ and unoriented $X$-homeomorphism $f\colon  N\rightarrow N'$, we have $\tau(W)=\tau(W_{1})\circ f_{\sharp}\circ\tau(W_{0})$, where $W$ is the unoriented $X$-cobordism obtained from $W_{1}$ and $W_{0}$ by gluing along $f$. 
\end{defn}
In the case where $\pi$ is trivial, $X$-niceness is equal to {\it niceness} introduced in \cite{turner-turaev:2006}. 
Then the following lemma is an easy consequence of Theorem~$\ref{turaev}$. 
\begin{lem}[\cite{turner-turaev:2006}]\label{nice}
$(1)$
Let $(W_{0}, M_{0}, N_{0})$ be an unoriented $X$-cobordism obtained from two unoriented $X$-cobordisms $(W'_{0}, M_{0}, N'_{0})$ and $(W''_{0}, M''_{0}, N_{0})$ by gluing along an unoriented $X$-homeomorphism $g\colon N'_{0}\rightarrow M''_{0}$. If $W'_{0}$ and $W''_{0}$ are $X$-nice, so is $W_{0}$. 
\par
$(2)$
Let $(W_{0}, M_{0}, N_{0})$ and $(W_{1}, M_{1}, N_{1})$ be oriented $X$-cobordisms and $f\colon N_{0}\rightarrow M_{1}$ be an orientation preserving $X$-homeomorphism. Then we have $\tau(W)=\tau(W_{1})\circ f_{\sharp}\circ\tau(W_{0})$, where $W$ is the oriented $X$-cobordism obtained from $W_{1}$ and $W_{0}$ by gluing along $f$. 
\end{lem}
\begin{prop}\label{gluing}
The six unoriented $X$-cobordisms depicted in Figure~$\ref{basic cobordism}$ are $X$-nice. 
Hence the pair $(A, \tau)$ constructed in Theorem~$\ref{mainthm}$ from an extended crossed $\pi$-algebra $(L, \eta, \varphi, \{\theta_{\alpha}\}_{\alpha\in\pi}, \Phi)$ satisfies the axioms of Definition~$\ref{HQFT}$. 
\end{prop}
\begin{proof}
Let $(W_{0}, M_{0}, N, g)$ and $(W_{1}, N', M_{1}, g')$ be two unoriented X-cobordisms and $f\colon N\rightarrow N'$ be an unoriented X-homeomorphism of unoriented $X$-cobordisms. Let $(W_{1}\cup_{f} W_{0}, g'\cup_{f}g)$ be an unoriented cobordism obtained from $W_{0}$ and $W_{1}$ by gluing along $f$
.  
Moreover we suppose that we have $(W_{1}, g')\cong (W_{1}^{{\rm or}}, g'|_{W_{1}^{{\rm or}}})\sharp(\mathbf{R}P^{2}, 
f_{\alpha_{1} })\sharp\dots\sharp(\mathbf{R}P^{2}, f_{\alpha_{n} })$. 
\par
(I) The case where $(W_{0}, M_{0}, N, g)$ is a cobordism depicted in Figure~$\ref{basic cobordism}$ ($3$). 
\par
In this case, we can choose orientations of $W_{1}$ and $W_{0}$ such that $f$ is an orientation preserving homeomorphism. By Lemma~$\ref{nice}$, we have the following equation: 
\begin{align*}
\tau(W_{1}\cup_{f} W_{0}, g'\cup_{f}g)
=&\psi_{\alpha_{1}, \cdots, \alpha _{n}}\circ \tau(W_{1}^{{\rm or}}\cup_{f} W_{0})\\
=&\psi_{\alpha_{1}, \cdots, \alpha _{n}}\circ \tau(W_{1}^{{\rm or}})\circ f_{\sharp}\circ \tau(W_{0})\\
=&\tau(W_{1}, g')\circ f_{\sharp} \circ\tau(W_{0}, g). 
\end{align*}
Hence $(W_{0}, M_{0}, N, g)$ is $X$-nice. 
\par
(II) The case where $(W_{0}, M_{0}, N, g)$ is a cobordism depicted in Figure~$\ref{basic cobordism}$ ($1$), ($4$) and ($5$). 
\par
In this case, we can use the same proof as (I). 
\par
(III) The case where $(W_{0}, M_{0}, N, g)$ is a cobordism depicted in Figure~$\ref{basic cobordism}$ ($6$). 
\par
Suppose that $(W_{0}, M_{0}, N, g)$ is given an unoriented $X$-homeomorphism 
$(W_{0}, g)\cong (W_{0}^{{\rm or}}, g|_{W_{0}^{{\rm or}}})\sharp(\mathbf{R}P^{2}, f_{\alpha })$. Then we have 
\begin{align*}
(W_{1}\cup_{f} W_{0}, g'\cup_{f}g)\cong (W_{1}^{{\rm or}}\cup_{f} W_{0}^{{\rm or}}, g'|_{W_{1}^{{\rm or}}}\cup_{f}g|_{W_{0}^{{\rm or}}})\sharp(\mathbf{R}P^{2}, f_{\alpha })\sharp(\mathbf{R}P^{2}, f_{\alpha_{1} })\sharp\dots\sharp(\mathbf{R}P^{2}, f_{\alpha_{n} }). 
\end{align*}
Furthermore we can give orientations of $W_{1}^{{\rm or}}$ and $W_{0}^{{\rm or}}$ such that $f$ preserves the orientations. Then we have the following equation:
\begin{align*}
\tau(W_{1}\cup_{f} W_{0}, g'\cup_{f}g)
&=\psi_{\alpha}\circ\psi_{\alpha_{1}, \cdots \alpha _{n}}\circ\tau(W_{1}^{{\rm or}}\cup_{f} W_{0}^{{\rm or}})\\
&=\psi_{\alpha}\circ\psi_{\alpha_{1}, \cdots \alpha _{n}} \circ \tau(W_{1}^{{\rm or}})\circ f_{\sharp}\circ \tau(W_{0}^{{\rm or}})\\
&=\psi_{\alpha_{1}, \cdots \alpha _{n}}\circ \tau(W_{1}^{{\rm or}})\circ f_{\sharp}\circ\psi_{\alpha}\circ \tau(W_{0}^{{\rm or}})\\
&=\tau(W_{1}, g')\circ f_{\sharp}\circ\tau(W_{0}, g).
\end{align*}
The third equality follows from Figure~$\ref{tousiki1}$. The equation in Figure~$\ref{tousiki1}$ follows from the fact that $f_{\sharp}=$id or $\Phi$ and that $\Phi(\theta_{\alpha})=\theta_{\alpha}$ for any $\theta_{\alpha}\in \pi$. 
Hence $(W_{0}, M_{0}, N, g)$ is $X$-nice. 
\begin{figure}[!h]
\begin{center}
\includegraphics[scale=0.51]{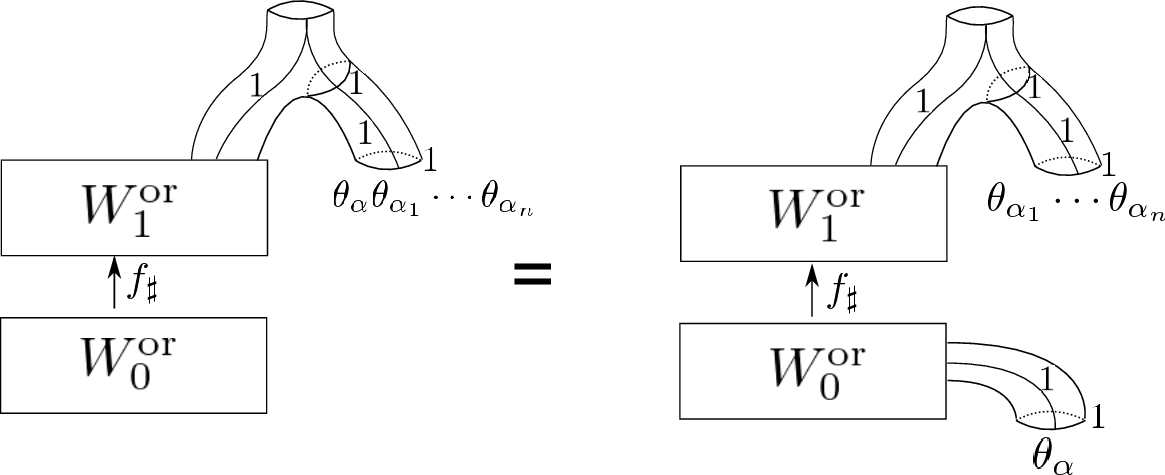}
\end{center}
\caption{$\psi_{\alpha}\circ\psi_{\alpha_{1}, \cdots \alpha _{n}} \circ \tau(W_{1}^{{\rm or}})\circ f_{\sharp}\circ \tau(W_{0}^{{\rm or}})
=\psi_{\alpha_{1}, \cdots \alpha _{n}}\circ \tau(W_{1}^{{\rm or}})\circ f_{\sharp}\circ\psi_{\alpha}\circ \tau(W_{0}^{{\rm or}})$. }
\label{tousiki1}
\end{figure}
\par
(IV) The case where $(W_{0}, M_{0}, N, g)$ is a cobordism depicted in Figure~$\ref{basic cobordism}$ ($2$). 
\par
If we can give orientations of $W_{1}^{{\rm or}}$ and $W_{0}$ so that $f$ preserves them, we can use the same argument as (I). 
Suppose that we cannot give such orientations. Then there are an unoriented $X$-cobordism $(W_{3}, M_{3}, N_{3})$, an unoriented $X$-cobordism $(W_{2},  M_{2}, N_{2})$ and an unoriented $X$-homeomorphism $f'\colon N_{2}\rightarrow M_{3}$ such that $W_{1}=W_{3}\cup _{f}W_{2}$ (see Figure~$\ref{W1}$), where $(W_{2}, M_{2}, N_{2})$ is unoriented $X$-homeomorphic to the unoriented $X$-cobordism depicted in Figure~$\ref{basic cobordism}$ (3).
Let $W_{4}$ be an unoriented $X$-cobordism $W_{2}\cup_{f}W_{0}$ (see Figure~$\ref{W4}$). 
\begin{figure}[!h]
\begin{center}
\includegraphics[scale=0.29]{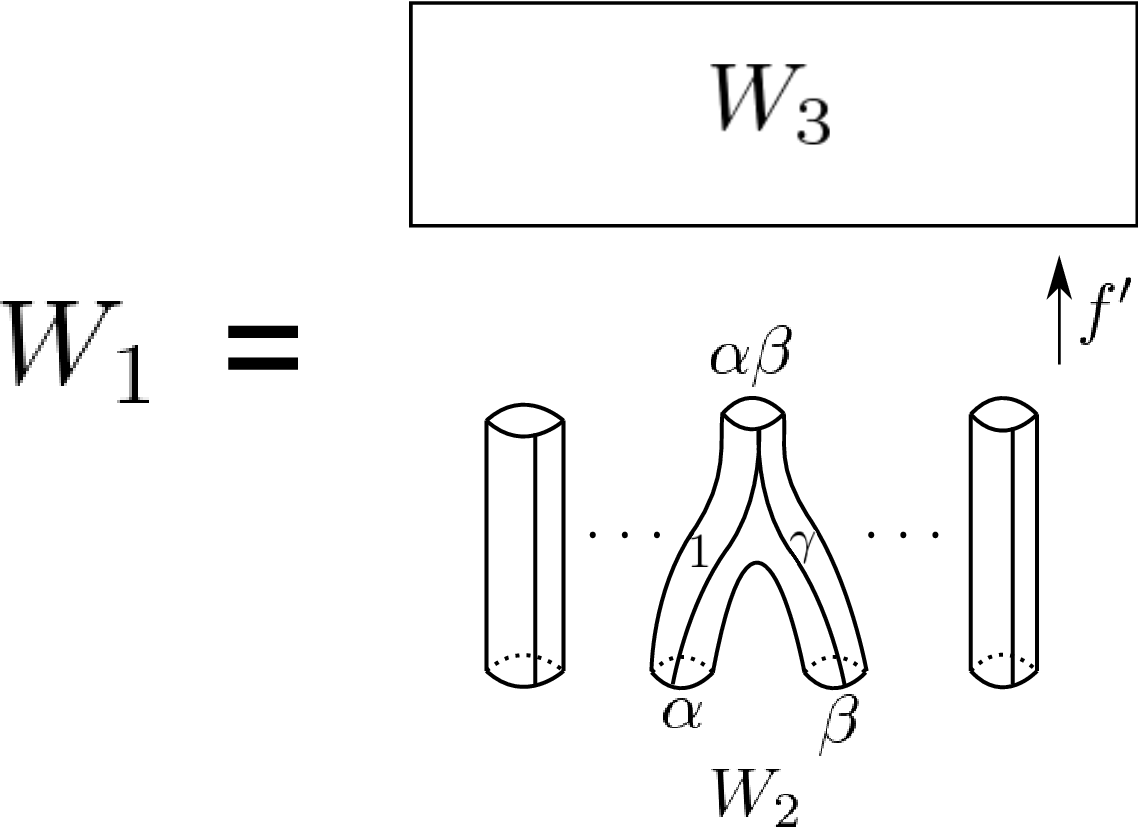}
\end{center}
\caption{$W_{1}=W_{3}\cup_{f'}W_{2}$. }
\label{W1}
\end{figure}
\begin{figure}[!h]
\begin{center}
\includegraphics[scale=0.29]{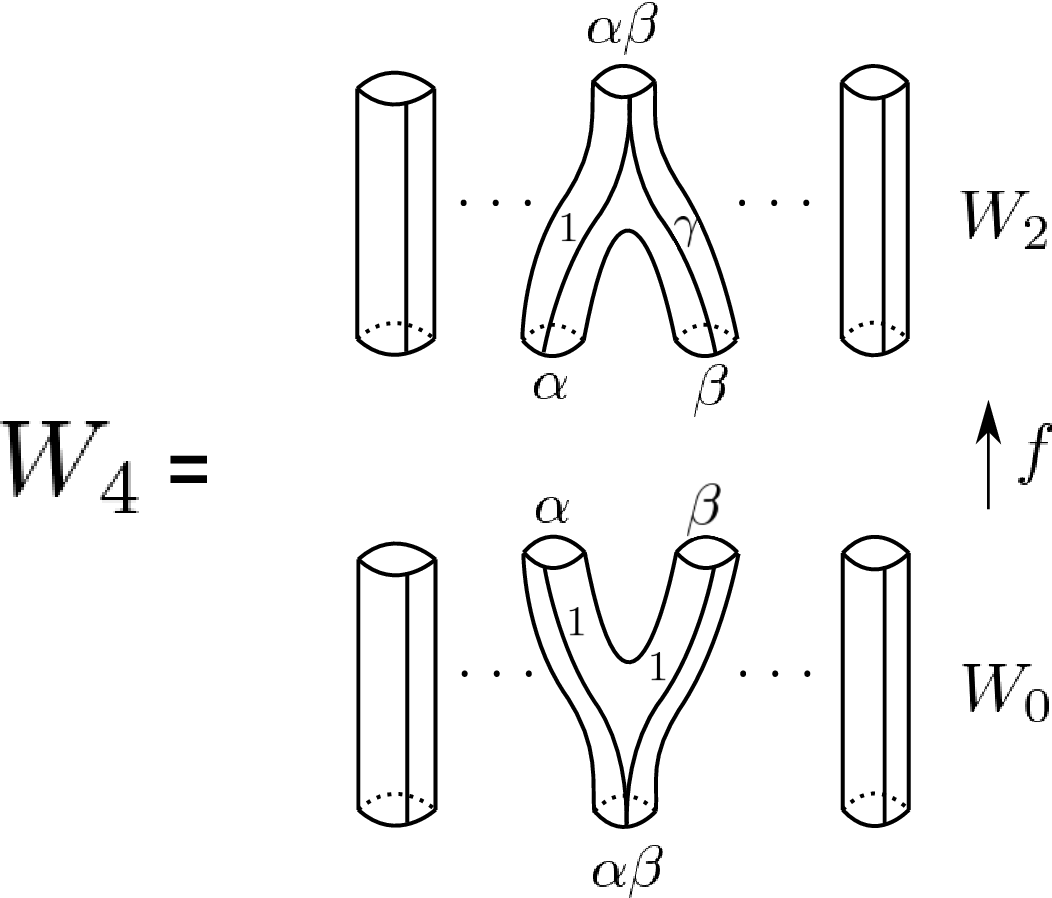}
\end{center}
\caption{$W_{4}=W_{2}\cup_{f}W_{0}$. }
\label{W4}
\end{figure}
It follows from the proof of Lemma~$\ref{lem8}$ that $(W_{4}, g)\cong (W_{4}^{{\rm or}}, g|_{W_{4}^{{\rm or}}})\sharp(\mathbf{R}P^{2}, f_{\alpha\gamma })\sharp(\mathbf{R}P^{2}, f_{\alpha})$. 
Moreover it follows from Lemma~$\ref{lem8}$ and the definition of extended crossed group-algebra (see Definition $\ref{extcross} (7)$) that 
\begin{center}
$\tau(W_{4})=\tau(W_{2})\circ f_{\sharp}\circ \tau(W_{0})$. 
\end{center}
Now we have 
\begin{multline*}
(W_{3}\cup_{f'} W_{4}, g'\cup_{f}g)\cong \\
(W_{3}^{{\rm or}}\cup_{f'} W_{4}^{{\rm or}}, (g'\cup_{f}g)|_{W_{3}^{{\rm or}}\cup_{f'} W_{4}^{{\rm or}}})\sharp(\mathbf{R}P^{2}, f_{\alpha\gamma })\sharp(\mathbf{R}P^{2}, f_{\alpha} )\sharp(\mathbf{R}P^{2}, f_{\alpha_{1} })\sharp\dots\sharp(\mathbf{R}P^{2}, f_{\alpha_{n} }) 
\end{multline*}
and we can give orientations of $W_{3}^{{\rm or}}$ and $W_{4}^{{\rm or}}$ so that $f'$ preserve them. Hence we have the following equation: 
\begin{align*}
\tau(W_{1}\cup_{f} W_{0}, g'\cup_{f}g)&=(W_{3}\cup_{f'} W_{4}, g'\cup_{f}g)\\
&=\psi_{\alpha_{1}, \cdots \alpha _{n}}\circ\psi_{\alpha\gamma}\circ\psi_{\alpha}\circ\tau(W_{3}^{{\rm or}}\cup_{f'} W_{4}^{{\rm or}})\\
&=\psi_{\alpha_{1}, \cdots \alpha _{n}}\circ\psi_{\alpha\gamma}\circ\psi_{\alpha}\circ\tau(W_{3}^{{\rm or}})\circ f'_{\sharp}\circ \tau(W_{4}^{{\rm or}})\\
&=\psi_{\alpha_{1}, \cdots \alpha _{n}}\circ\tau(W_{3}^{{\rm or}})\circ f'_{\sharp}\circ\psi_{\alpha\gamma}\circ\psi_{\alpha}\circ\tau(W_{4}^{{\rm or}})\\
&=\tau(W_{3})\circ f'_{\sharp}\circ \tau(W_{4})\\
&=\tau(W_{3})\circ f'_{\sharp}\circ\tau(W_{2})\circ f_{\sharp}\circ \tau(W_{0})\\
&=\tau(W_{1})\circ f_{\sharp}\circ \tau(W_{0}). 
\end{align*}
Hence $W_{0}$ is $X$-nice. 
\end{proof}
\section{Examples}\label{ex}
In this section, we construct examples of HQFTs and extended crossed group algebras. 
\par
Firstly we will construct an example of unoriented HQFTs.
\begin{example}
 This construction is similar to ``primitive cohomological HQFT" constructed by Turaev \cite{turaev:1999} and his construction is inspired by the work of Freed and Quinn \cite{F-Q}. 
Let $\pi$ be $\mathbf{Z}/2\mathbf{Z}$ and $X$ be a $K(\pi, 1)$ space (in particular $X$ is homotopy equivalent to $\mathbf{R}P^{\infty}$). 
Given $d\geq 0 $ we take a $(d+1)$-dimensional cocycle $\theta \in C^{d+1}(X; R^{\times})$, where $R^{\times}$ is the unit group of $R$. 
For any unoriented $d$-dimensional $X$-manifold $(M, g)$, we define an $R$-module $A(M,g)$
 by $Rv_{a}$, where $a\in C_{d}(M; \mathbf{Z}/2\mathbf{Z})$ is a fundamental cycle and $Rv_{a}$ is the free $R$-module of rank $1$ generated by $v_{a}$. 
If $a, b\in C_{d}(M; \mathbf{Z}/2\mathbf{Z})$ are two fundamental cycles, then we give the relation $v_{a}=g^{\ast}(\theta )(c)v_{b}$, where $c$ is a $(d+1)$-dimensional
singular chain in $M$ such that $\partial{c}=a+b$. The element $g^{\ast}(\theta)(c)\in R^{\times}$ does not depend on the choice of $c$. 
For any unoriented $X$-homeomorphism $f\colon (M,g)\rightarrow (M',g')$, we define an $R$-homomorphism $f_{\sharp}\colon A(M,g)\rightarrow A(M',g')$ by $f_{\sharp}(v_{a})=v_{f_{\ast}(a)}$. 
\par
Let $(W, M_{0}, M_{1}, g)$ be an unoriented $(d+1)$-dimensional $X$-cobordism. 
Take a cycle $B\in C_{d+1}(W,\partial{W}; \mathbf{Z}/2\mathbf{Z})$ such that $[B]\in H_{d+1}(W, \partial{W}; \mathbf{Z}/2\mathbf{Z})$ is the fundamental class. 
Then we have $\partial{B}=a_{0}+a_{1}$, where $\partial\colon  C_{d+1}(W, \partial{W}; \mathbf{Z}/2\mathbf{Z})\rightarrow C_{d}(M_{0}; \mathbf{Z}/2\mathbf{Z})\oplus C_{d}(M_{1};\mathbf{Z}/2\mathbf{Z})$ is the connected homomorphism and $a_{0}\in C_{d}(M_{0}; \mathbf{Z}/2\mathbf{Z})$, $a_{1}\in C_{d}(M_{1};\mathbf{Z}/2\mathbf{Z})$ are fundamental cycles. 
Then we define an $R$-homomorphism $\tau(W, g)\colon A(M_{0}, g_{M_{0}})\rightarrow A(M_{1}, g_{M_{1}})$ by $\tau(W, g)(v_{a_{0}})=(g^{\ast}(\theta )(B))^{-1}v_{a_{1}}$. 
The map $\tau(W, g)$ does not depend on the choice of $B$. 
\par
The pair $(A, \tau)$ is an unoriented $(d+1)$-dimensional HQFT with target $X$. 
Moreover the isomorphism class $(A, \tau)$ does not depend on the choice of a singular cocycle representation $\theta $ of the homology class $[\theta ]\in H^{d+1}(X,;R^{\times})$. 
For any closed unoriented $(d+1)$-dimensional $X$-cobordism $(W, g)$, the map $\tau(W, g)$ is an involution. 
\end{example}
Secondly we make an example of extended crossed group algebras below. 
\begin{example}
Let $\pi$ be the group $\mathbf{Z}/2\mathbf{Z}=\{1, -1\}$ and $\{l_{\alpha}\}_{\alpha\in\pi}$ be a set whose index set is $\pi$. 
Let $\{\kappa_{\alpha , \beta }\in R^{\times}\}_{\alpha , \beta \in\pi}$ be a normalized $2$-cocycle, that is, $\kappa_{1,1}=1$ and $\kappa_{\alpha , \beta }\kappa_{\alpha \beta , \gamma }=\kappa_{\alpha , \beta \gamma }\kappa_{\beta ,\gamma }$, where $R^{\times}$ is the group of units of $R$. 
Note that for any $\alpha \in \pi$, we have $\kappa_{1, \alpha }=\kappa_{\alpha, 1 }=1$. 
\par
For any $\alpha\in\pi$, let $L_{\alpha}$ be the free $R$-module of rank $1$ generated by $l_{\alpha}$, that is, $L_{\alpha}=Rl_{\alpha}$. 
Put $L=L_{1}\oplus L_{-1}$. 
Multiplication of $L$ is defined by $l_{\alpha}l_{\beta}=\kappa_{\alpha , \beta }l_{\alpha\beta}$. 
A bilinear form $\eta \colon L\otimes L\rightarrow R$ is defined by $\eta(l_{\alpha}\otimes l_{\alpha})=\kappa_{\alpha , \alpha }$ for any $\alpha\in\pi$ and $\eta(l_{\alpha}\otimes l_{\beta})=0$ for $\beta\neq\alpha$. 
For any $\beta \in\pi$, put $\varphi _{\beta }=\id$. 
Take an element $a\in R$ which satisfies $a^{2}=1$ and put $\theta_{\alpha}=al_{1}$ for any $\alpha\in\pi$. 
Then $(L=\bigoplus L_{\alpha}, \eta, \varphi, \{\theta_{\alpha}\}_{\alpha\in\pi}, \Phi)$ is an extended crossed $\pi$-algebra. 

We can easily prove that the algebra $(L=\bigoplus L_{\alpha}, \eta, \varphi, \{\theta_{\alpha}\}_{\alpha\in\pi}, \Phi)$ satisfies the axioms  in Definition~$\ref{extcross}$ except $(3), (7)$ and $(11)$. 
Since we have $l_{\alpha }l_{\beta }=\kappa_{\alpha ,\beta }l_{\alpha \beta }=\kappa_{\alpha ,\beta }l_{\beta \alpha }=\kappa_{\alpha ,\beta }\kappa_{\beta, \alpha }^{-1}l_{\beta }l_{\alpha }=l_{\beta }l_{\alpha }$, $L$ satisfies the axiom $(3)$. 
\par
To check the axiom $(7)$, we need to compute $\Delta _{\alpha, \beta}(l_{\alpha\beta})$ for any $\alpha, \beta\in\pi$. 
Put $\Delta _{\alpha, \beta}(l_{\alpha\beta})=kl_{\alpha}\otimes l_{\beta}$. 
Then we have 
\begin{align}
(\id\otimes\eta)\circ(\Delta _{\alpha , \beta}\otimes\id)(l_{\alpha \beta }\otimes l_{\beta })=l_{\alpha \beta }l_{\beta } \label{eq6-1}
\end{align}
{\rm(}see Figure~$\ref{computeyoseki}${\rm)}. 
The left hand side of ($\ref{eq6-1}$) is equal to $k\kappa_{\beta , \beta }l_{\alpha}$ and the right hand side is equal to $\kappa _{\alpha \beta , \beta }l_{\alpha\beta \beta }=\kappa _{\alpha \beta , \beta }l_{\alpha}$. Hence $k=\kappa _{\beta, \beta}^{-1}\kappa _{\alpha\beta, \beta}$ and we have 
\begin{align*}
m\circ(\Phi \otimes \varphi _{\gamma})\circ\Delta _{\alpha, \beta}(l_{\alpha\beta})&=l_{\alpha \beta }. 
\end{align*}
\par
To check the axiom $(11)$, we need to compute $q(1)\in L_{1}$. 
We consider $\tau^{L}(Q')$, where the cobordism $Q'$ is depicted in Figure~$\ref{decomQ}$ whose bottom base is empty and whose top base is $(\mathbf{S}^{1}, \alpha \beta )\sqcup(\mathbf{S}^{1}, \alpha \beta )$. 
Put $\tau^{L}(Q')(1)=k'l_{\alpha \beta }\otimes l_{\alpha \beta }$. 
Now we have 
\begin{align}
(\id\otimes\eta)\circ(\tau^{L}(Q')\otimes \id)(l_{\beta \alpha })=\varphi _{\beta \gamma }(l_{\beta \alpha }) \label{eq6-2}
\end{align}
{\rm(}see Figure~$\ref{comp2}${\rm)}.
The left hand side of ($\ref{eq6-2}$) is equal to $k'\kappa _{\alpha \beta , \alpha \beta }l_{\alpha\beta }$ and the right hand side is equal to $l_{\alpha\beta }$. 
Hence $k'=\kappa _{\alpha \beta, \alpha \beta}^{-1}$ and we have $q(1)=m\circ\tau^{L}(Q')(1)=k'l_{\alpha \beta }l_{\alpha \beta }=l_{1}$. 
This shows that $L$ satisfies the axiom $(11)$. 
\begin{rem}
Note that Turaev \cite{turaev:1999} shows that the algebra $(L=\bigoplus L_{\alpha}, \eta, \varphi)$ is a crossed $\pi$-algebra. 
\end{rem}
\begin{figure}[!h]
\begin{center}
\includegraphics[scale=0.25]{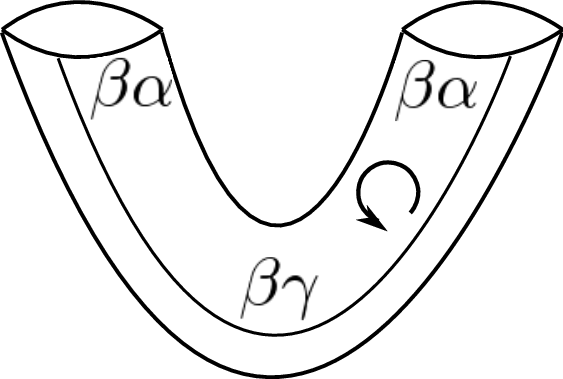}
\end{center}
\caption{The cobordism $(Q', \emptyset, (\mathbf{S}^{1}, \alpha \beta )\sqcup(\mathbf{S}^{1}, \alpha \beta ))$.}
\label{decomQ}
\end{figure}
\end{example}
\par
\noindent
{\bf Acknowledgement: } The author would like to express his sincere gratitude to Hitoshi Murakami for his encouragement. 
\bibliographystyle{amsplain}
\bibliography{mrabbrev,tagami}
\end{document}